% FILE intro.tex
% put the files intro.tex makro.tex litp.tex 
% in the same directory and execute (twice) the command: tex intro
\magnification=1200
% FILE makro.tex
\newif\iftestmodus  % testmodus: Nur einfache Ausgabe der Labels
\newif\iflabelanzeige % Sorgt daf"ur, da"s bei ausgeschaltetem Testmodus die
                    % Labels angezeigt werden (bei den S"atzen, nicht bei
                    % Referenzierungen).
\newif\ifkapkopfzeile % Anzeige des Kapitels in der Kopfzeile
\newif\ifseitenzahl % Anzeige der Seitenzahl in der Kopfzeile
\newif\ifdatum % Anzeige des Datums in der Kopfzeile (statt des Kapitels)
\newif\ifdoppel % abwechselnd links und rechts Seitennummer

% Voreinstellungen
\hfuzz=1pt
\newdimen\anteil
\anteil=0pt
\voffset=1.5cm
\topskip=8mm
\frenchspacing
\testmodusfalse
\labelanzeigetrue
\kapkopfzeilefalse
\datumfalse
\doppelfalse
\seitenzahltrue
%Zeilenabstand
\advance\baselineskip by 0.9mm

% gr"o"sere Schriften
\font\rmeins=cmr10 scaled\magstep1
\font\rmzwei=cmr10 scaled\magstep2
\font\rmdrei=cmr10 scaled\magstep3
\font\ieins=cmmi10 scaled\magstep1
\font\izwei=cmmi10 scaled\magstep2
\font\idrei=cmmi10 scaled\magstep3
\font\syeins=cmsy10 scaled\magstep1
\font\syzwei=cmsy10 scaled\magstep2
\font\sydrei=cmsy10 scaled\magstep3
\font\exdrei=cmex10 scaled\magstep3
\font\exzwei=cmex10 scaled\magstep2
\font\exeins=cmex10 scaled\magstep1
\def\Roman#1{\uppercase\expandafter{\romannumeral#1}}

\font\zusatz=msbm10 scaled 1728
\font\zusatzz=msbm10 scaled 1440

% Querverweise

% Einsch"ube von LaTeX: Arbeitet mit beliebigen Zeichenketten.
\def\namedef#1{\catcode`\@=11 \expandafter\def\csname #1\endcsname\catcode`\@=12 }% entspricht \def
\def\nameuse#1{\catcode`\@=11 \csname #1\endcsname\catcode`\@=12 }% Ausf"uhren des mit \namedef Definierten
\newtoks\EEEzk  % z. B. "Proposition"
\newcount\EEEa  % Kapitelnummer
\newcount\EEEb  % Unterkapitelnummer
\newcount\EEEc  % Nummer des Satzes usw.
\newcount\EEEs % Seitennummer
\newcount\Csection % aktuelle Kapitelnummer
\newcount\Csubsection % aktuelle Unterkapitelnummer
\newcount\Csatz % Satznummer
\newcount\Cglei % Gleichungsnummer
%
% Abschlu"sbefehl:

%
% Wirkungsweise der Labels: Setzen von Variablen
\def\Definieren#1#2#3#4#5#6{
\EEEa=-3000 \nameuse{r@#1}\relax%
\ifnum\EEEa=-3000\else\message{!!! Label '#1' beim letzten Durchlauf doppelt.}\fi
\namedef{r@#1}{%
 \Setzen{#2}{#3}{#4}{#6}{\ignorespaces{#5}}}}
% Setzbefehl:
\def\Setzen#1#2#3#4#5{\global\EEEzk={#5}\global\EEEa=#1%
\global\EEEb=#2\global\EEEc=#3\global\EEEs=#4}
%
% Label:
%
\def\Label#1#2#3#4#5% #1: Labelname #2: Kapitel, #3: Unterkapitel, #4: Nummer
% #5: Bezeichnung wie "Proposition"
{\edef\Schreiben{\write\Meine{\string\Definieren{#1}{\number#2}{\number#3}{\number#4}{#5}
{\noexpand\folio}}}\Schreiben}
% d. h. in Datei schreiben
%
%
% Label in Gleichung:
\def\Eqlabel#1#2{\iftestmodus(#1)\else%
\global\advance\Cglei by 1
\Label{#1}\Csection\Csubsection\Cglei{Gleichung} % Verzeichniseintrag
#2 (\number\Cglei) %Gleichungsnummer anschreiben
\iflabelanzeige\rlap{[#1]}\fi\fi}

% mit \eqno-Befehl:
\def\Eqno#1{\Eqlabel{#1}\eqno}
% ohne (f"ur \eqalignno}:

%
% Label in Satz:
\def\SLabel#1#2%  #2: "Proposition" #1: Labelname
{\iftestmodus$\{#1\}$\else%
\global\advance\Csatz by 1
\Label{#1}\Csection\Csubsection\Csatz{#2}% Verzeichniseintrag
\bf\iflabelanzeige\textindent{[#1]}\fi%Fett, evtl. Labelanzeige
\number\Csubsection.\number\Csatz.\fi% Satznummer aufschreiben
\ #2\ifdoppelpunkt:\else\doppelpunkttrue\fi}% "Proposition", Doppelpunktsetzung
%
% Referierung: Mit Ausgabe von "Proposition" usw.:
%
\newif\iffehler

% \Refgrundlage setzt die Variablen auf die richtigen Werte bzw.
% setzt \fehlertrue, falls der Name nicht definiert ist.
\def\Refgrundlage#1{\global\EEEa=-3000\nameuse{r@#1}\ifnum\EEEa=-3000%
\global\fehlertrue\message{!!! Label '#1' auf \folio. Seite nicht definiert.}%
\else\global\fehlerfalse\fi}

% Anwendungen von \Refgrundlage:

% Seitenausgabe
\def\Seitenref#1{\iftestmodus[#1]\else%
\Refgrundlage{#1}%
\iffehler !![#1]\else%
\number\EEEs\fi\fi}

% Ausgabe in interner Numerierung

\def\Ref#1{\iftestmodus[#1]\else%
\Refgrundlage{#1}%
\iffehler !![#1]\else\the\EEEzk~\ifnum\EEEa=\Csection\else\uppercase%
\expandafter{\romannumeral\EEEa}.\fi\number\EEEb.%
\number\EEEc\fi\fi}
%
% F"ur Gleichungen:
\def\Refg#1{\iftestmodus Gleichung [#1]\else%
\Refgrundlage{#1}%
\iffehler !![#1]\else\the\EEEzk~(\number\EEEc)\fi\fi}

\let\SRefg=\Srefg
%
% mit l:
\def\Refgl#1{\iftestmodus Gleichung [#1]\else%
\Refgrundlage{#1}%
\iffehler !![#1]\else\the\EEEzk~$(\number\EEEc)_l$\fi\fi}

% F"ur Kapitel:
\def\Refk#1{\iftestmodus[#1]\else%
\Refgrundlage{#1}%
\iffehler !![#1]\else\the\EEEzk~\Roman\EEEa\ifnum\EEEb>-3000\relax
\number\EEEb\fi\fi\fi}
\def\Srefk#1{\iftestmodus[#1]\else%
\Refgrundlage{#1}%
\iffehler !![#1]\else\Roman\EEEa\ifnum\EEEb>-3000\relax
\ifnum\EEEa>0\relax.\fi\number\EEEb\fi\fi\fi}
% Ohne Ausgabe:
%
\def\Sref#1{\iftestmodus[#1]\else%
\Refgrundlage{#1}\iffehler !![#1]\else\ifnum\EEEa=\Csection
\else\uppercase\expandafter{\romannumeral\EEEa.}\fi\number\EEEb.\number\EEEc\fi\fi}
\let\SRef=\Sref
%%
% Anfangsbefehl:
%
\newread\Meine
\immediate\openin\Meine=\jobname.hlf
\newif\ifnichtfertig  % Datei zu Ende?
\nichtfertigtrue % noch nicht
% Datei bis zum Ende einlesen
\loop \ifeof\Meine\nichtfertigfalse\fi
\ifnichtfertig
\immediate\read\Meine to \EEEV
\EEEV
\repeat
\immediate\closein\Meine
% zum Schreiben wieder "offnen:
\newwrite\Meine
\immediate\openout\Meine=\jobname.hlf
%

% ANDERES

% Seiten"uberschrift
\newtoks\Ueberschr
\newtoks\LUeberschr
\Ueberschr={}
\LUeberschr={}
\newcount\stunden
\newcount\minuten
\stunden=\time \divide\stunden by 60
\minuten=-\stunden
\multiply \minuten by 60
\def\datumtageszeit{{\number\day. \number\month. \number\year,
\number\stunden :\ifnum\minuten<10\relax 0\fi\number\minuten}}
\advance\minuten by \time
\nopagenumbers
\headline={\tenrm %Datei "\jobname"
\ifkapkopfzeile\ifdoppel\ifodd\pageno\the\Ueberschr\hss\ifseitenzahl\quad\folio\fi
\else\ifseitenzahl\folio\fi\quad\hss\the\LUeberschr\fi\else
\the\Ueberschr\hss\ifseitenzahl\quad\folio\fi\fi
\else\hfil\global\kapkopfzeiletrue\fi
\global\seitenzahltrue}

% deutsche Anf"uhrungszeichen

     % !!!!!!!!!!!!!!!! Scheinparameter!

% Verwaltung der S"atze

\newif\ifsatz
\satzfalse
\newif\ifdoppelpunkt
\doppelpunkttrue

% Vorkehrungen bei Satzanfang, Satzende
\def\Satzanfang{\removelastskip\goodbreak\medskip\par \anteil=2pt\parskip=0pt
\parindent=0pt
% \global\leftskip=10pt \hglue -10pt
 }

\def\Satzende{\anteil=0pt\parskip=2pt \global\leftskip=0pt\parindent=10pt\goodbreak\medskip\rm}

\def\Einl#1#2{\ifsatz\message{!! Satz vor #1 nicht beendet!}
\else\satztrue\fi\Satzanfang\SLabel{#1}{#2}}

\def\Def#1{\Einl{#1}{Definition}\rm}
\def\Lemma#1{\Einl{#1}{Lemma}\rm\sl}

\def\DefLemma#1{\Einl{#1}{Definition-Lemma}\rm\sl}

\def\Prop#1{\Einl{#1}{Proposition}\rm\sl}

\newif\ifbemerkung
\def\Bew{\ifbemerkung\bemerkungfalse\else
\removelastskip\smallbreak\fi
\ifsatz
\rm Beweis\ifdoppelpunkt:\else\doppelpunkttrue\fi\ %
\else\message{Beweis in falscher Umgebung!}\fi}

\def\Ende{\ifsatz\satzfalse\Satzende\rm\else\message{!! Falsches Ende}\fi}

% Latexbefehle/noch zu definieren. Parameter: soll als Referenz gegeben
% werden, und falls bestimmtes "Probesymbol" definiert ist, soll
% Parameter mit ausgegeben werden

\def\Anzeige#1{\ifdatum\Ueberschr=\datumtageszeit\else\Ueberschr={Chapter \Roman\Csection\ifnum\Csection>0.\fi\number\Csubsection. #1}\fi}
\def\Seitenende{\vfill\eject}
\def\Umbruchbefehl{\relax}
\def\Umbruch{\Umbruchbefehl
 \if\Umbruchbefehl\relax\def\Umbruchbefehl{\Seitenende}\fi}
\def\section#1#2{\Umbruch\def\Umbruchbefehl{\relax}\global\kapkopfzeilefalse%
 \noindent%
 {\iftestmodus\else\global\advance\Csection by 1%
 \global\Csubsection=0%
 \Label{#1}\Csection{-3000}{-3000}{Teil}%
 \rmdrei\everymath{\textfont\msbfam=\zusatz\textfont0=\rmdrei\scriptfont0=\rmzwei
 \scriptscriptfont0=\rmeins\textfont1=\idrei\scriptfont1=\izwei
 \scriptscriptfont1=\ieins\textfont2=\sydrei\scriptfont2=\syzwei
 \scriptscriptfont2=\syeins\textfont3=\exdrei\scriptfont3=\exzwei
 \scriptscriptfont3=\exeins}%
\uppercase\expandafter{\romannumeral\Csection}.  \fi #2} \rm \par}
\def\subsection#1#2{\Umbruch\bigbreak\Anzeige{#2}\global\kapkopfzeilefalse\noindent\seitenzahlfalse%
{\rmzwei%
 \iftestmodus\else\global\advance\Csubsection by 1
 \global\Csatz=0 \number\Csubsection.
 \Label{#1}\Csection\Csubsection{-3000}{Kapitel}\fi%
\everymath{\textfont\msbfam=\zusatzz
 \textfont0=\rmzwei\scriptfont0=\rmeins
 \scriptscriptfont0=\tenrm\textfont1=\izwei\scriptfont1=\ieins
 \scriptscriptfont1=\teni\textfont2=\syzwei\scriptfont2=\syeins
 \scriptscriptfont2=\tensy\textfont3=\exzwei\scriptfont3=\exeins
 \scriptscriptfont3=\tenex}#2 }\rm \par}
\def\Abschnitt#1{\removelastskip\bigbreak\rmeins
{\everymath{\textfont\msbfam=\zusatzz
 \textfont0=\rmeins\scriptfont0=\sevenrm
 \scriptscriptfont0=\tenrm\textfont1=\ieins\scriptfont1=\ieins
 \scriptscriptfont1=\teni\textfont2=\syeins\scriptfont2=\syeins
 \scriptscriptfont2=\tensy\textfont3=\exeins\scriptfont3=\exeins
 \scriptscriptfont3=\tenex}\noindent #1 }\rm \par}

% Definition von ANH: zum Platzausgleich: hat nachfolgendes Zeichen den Wert
% eines Satzzeichens - .,?!-~;: - oder "Klammer zu" oder "\ " (d. h.
% erzwungener Leerraum), dann wird kein
% Leerraum gelassen und der Parameter wieder eingesetzt. Ansonsten
% wird ein Leerraum gelassen und der Parameter eingesetzt. Dies
% sorgt f"ur richtigen Leerraum nach dem Makro. Verwendung von ANH:
% "ublicherweise ohne Parameter (nach Makroexpansion schnappt sich
% ANH das n"achste Zeichen automatisch).
% soll nach ANH entgegen dieser Regeln ein Leerraum kommen, so setzt
% man zun"achst einen erzwungenen Leerraum dahinter: "\ ",
% soll kein Leerraum kommen, so setzt man den Befehl "\empty" dahinter.
\def\ANH{\def\qpzy{\ifx\qpyz\empty\else\ifx,\qpyz\else\ifx\qpyz.\else\ifx\qpyz:\else\ifx\qpyz;%
\else\ifx\qpyz!\else\ifx\qpyz?\else\ifx\qpyz)\else\ifx\qpyz-\else%
\ifx\qpyz\space\else\ifx\qpyz~\else\space%
\fi\fi\fi\fi\fi\fi\fi\fi\fi\fi\fi}\futurelet\qpyz\qpzy}
\def\ANHB{\def\qpzy{\ifx\qpyz\empty\else\ifx\qpyz,\else\ifx\qpyz.\else\ifx\qpyz:\else\ifx\qpyz;%
\else\ifx\qpyz!\else\ifx\qpyz?\else\ifx\qpyz)\else\ifx\qpyz-\else%
\ifx\qpyz\ \else\ifx\qpyz~\else~%
\fi\fi\fi\fi\fi\fi\fi\fi\fi\fi\fi}\futurelet\qpyz\qpzy}

% Befehle mit S c h e i n parameter (f"ur richtigen Leerraum)

% Wurzel aus q

% hbox und mathematischer Modus

% Rechtscoidealunteralgebra

% Rechtscoidealunteralgebren

% Stern-Hopfalgebra

% M-Zahl (\alpha\beta-Zahl)

% N-Zahl (\alpha\gamma-Zahl)

% M- und N-Zahl

% Befehle  o h n e   S c h e i n parameter
% nur mathematischer Modus
 % Tensorzeichen
\let\tensor=\otimes
 % Epsilon
\let\eps=\varepsilon
 % Wurzel aus q

 % \alpha, \beta, \gamma, \delta

 % k ohne Null

 % A_{m,n}

 % <X>: Diamantenlemma

 % Unterstrich/Bogen    %   !!!!!   vielleicht noch umdefinieren.

 % schr"age Punkte
\def\dddots{\mathinner{\mkern1mu\raise7pt\vbox{\kern7pt\hbox{.}}\mkern2mu
    \raise5.5pt\hbox{.}\mkern2mu\raise4pt\hbox{.}\mkern2mu
    \raise2.5pt\hbox{.}\mkern2mu\raise1pt\hbox{.}\mkern1mu}}

% Zahlbereiche (mit AMSTEX), hiervor " als Buchstaben setzen:
\newfam\msbfam
\font\amsten=msbm10
\font\amsseven=msbm7
\font\amsfive=msbm5
\textfont\msbfam=\amsten
\scriptfont\msbfam=\amsseven
\scriptscriptfont\msbfam=\amsfive
\def\Bbb#1{{\fam\msbfam #1}}
\def\Zahlbereich#1{\ifmmode{\Bbb #1}\else$\Bbb #1$\fi}
 % {\cal R)

 % ganze Zahlen

% nat"urliche Zahlen
\def\N{\Zahlbereich{N}\ANH}
% komplexe Zahlen
\def\C{\Zahlbereich{C}\ANH}

 % q hoch -1
\def\qm{q^{-1}}
 % Bilinearform

 % Coinvariante Elemente

 % Kleiner-Gleich in Mengen X und Y

 % Folgepfeil, Pfeil
\def\Folgt{\ \Rightarrow \ }

 % Kern, Identit\"at, Hom
\def\Ke{{\rm Ke}}
\def\id{{\rm id}}
\def\Hom{{\rm Hom}}
 % Schnitt

 % Element mit eingeklammertem Index (zur Unterscheidung zu Sigmanotation)

 % 1/2
\def\halb{{1 \over 2}}
\def\einhalb{\mathchoice{{\textstyle\halb}}{\halb}{\halb}{\halb}}
 % Quadratwurzeln:

 % q-Binomialkoeffizient (ohne Index)

 % ... mit Index q^2, q^{-2}:

 % ... mit Index \lambda:

 % ... mit Wurzel

 % Terme in V4.9:

\def\vwover#1{{\if#11\else#1\fi vw\over (q-\qm)^2}}

 % Matrizen mit a b c d:

 % Bilinearformwerte

 % t_ij^l Matrixelemente

\def\a#1 #2 {a_{#1#2}}
 % in V4.6: B_{X'}

 % gro{\ss}e Matrix

 % \tilde A .. \tilde D

% Pfeile
\def\llongrightarrow{\relbar\joinrel\relbar\joinrel\longrightarrow}

% Pfeile mit "Verzierung"

% \"Ubereinanderstehendes

% Abk\"urzungen f\"ur Harpunenpfeile

% Textbefehle ohne Scheinparameter (TEXTMODUS)
 % Implikationen, enthalten (f\"ur Beweisgliederung)

% Trennungen
\hyphenation{Al-ge-bren-ho-mo-mor-phis-mus}
\hyphenation{Hopf-al-ge-bra}
\hyphenation{Hopf-al-ge-bren}
\hyphenation{Hopf-al-ge-bren-ho-mo-mor-phis-mus}

% Items, Abs\"atze (umdefiniert)
\parindent=0pt
\advance\parskip by 2pt

\def\voritem#1 {\hangafter=1\hangindent20pt{\rm #1 }}
\def\voritemitem#1 {\hangafter=1\hangindent30pt{\rm #1 }}
\def\item#1 {\vskip\anteil\par \hangafter=0 \hangindent20pt \textindent{#1}}
\def\itemitem#1 {\vskip\anteil\par \hangafter=0 \hangindent30pt \textindent{#1}}

\def\textindent#1{\indent\llap{{\rm #1}\enspace}\ignorespaces}

% references
\def\wTmP#1#2#3{[{#2}\def\test{#3}\ifx\test\empty\else, #3\fi]%
\expandafter\gdef\csname #1\endcsname{\wtMP{#2}}}
\def\wtMP#1#2{[{#1}\def\test{#2}\ifx\test\empty\else, #2\fi]}
\def\Vorber#1{{\def\lit##1 ##2 ##3{\expandafter
\gdef\csname ##1\endcsname{\wTmP{##1}{##2}}}
\gdef\wTMp{#1}\input #1 }}
\def\Litverz{{\def\wTmP##1##2##3{}
\def\lit##1 ##2 ##3{\csname ##1\endcsname{##3\iflabelanzeige{(##1)}\fi}{}{}}
\def\wtMP##1##2{\par\hskip0pt\hangindent=1.3cm\hangafter=1\hbox to 1cm{[##1]\hfill} ##2}
\def\ww##1{\input ##1 }\ww{\wTMp} }}

\font\eufm=eufm10
\def\frak#1{\hbox{\eufm #1}}
\hsize=6.5true in
\vsize=9.2true in
\topskip=6 true mm
\parindent=10pt
\advance\voffset by -3 true mm
\advance\hoffset by -3 true mm
%\advance\leftskip by -3 true mm
%\advance\rightskip by 3 true mm 

\datumfalse
\labelanzeigefalse
\doppelfalse
\LUeberschr={Introduction to Quantum Groups}

\def\qgs{quantum groups\ANH}
\def\Quote#1{``#1''}
\def\R{\Zahlbereich{R}\ANH}
\let\rho=\varrho
% FILE litp.tex
\def\Quote#1{{\sl #1}}
\lit Abe A {E. Abe, \Quote{Hopf algebras,} Cambridge University Press 1980}
\lit Dix D {J. Dixmier, \Quote{$C^*$-Algebras,} North Holland 1977}
\lit Drin Dr {V.G. Drinfeld, \Quote{Quantum groups,} Proceedings ICM-1986,
Berkeley, AMS, 1987, pp. 798--820}
\lit DV DV {M. Dubois-Violette, G. Launer, The quantum group of a non-degenerate bilinear form, Phys. Lett. B {\bf 245} (1990), 175-177}
\lit Hay H {T. Hayashi, \Quote{Quantum deformation of classical groups,} 
Publ. RIMS, Kyoto {\bf 28} (1992),  57--81}
\lit KP KP {P. Kondratowicz, P. Podle\'s, \Quote{On representation theory
of quantum $SL_q(2)$ groups at roots of unity}, hep-th/9405079, 
   \Quote{Quantum groups and quantum spaces},
   R. Budzy\'nski et al. (eds.),
    Banach Center Publications {\bf40} (1997), 223--248,
    Inst. of Math., Polish Acad. Sci.}
%in {\sl Proceedings of the Banach minisemester on quantum groups and quantum
%spaces,} Warsaw 1995, Banach Center Publications; cf. also S. L.~%
%Woronowicz, a course \Quote{Quantum groups} at Faculty of
%Physics, University of Warsaw (1990/91)} 
\lit Koo Ko {T.H. Koornwinder, \Quote{General compact quantum groups, a tutorial,}
hep-th/9401114, in: \Quote{Representations of Lie groups and
quantum groups}, V.~Baldoni \& M.~A.~Picardello (eds.), Pitman
Res. Notes Math. Ser. {\bf311}, Longman Scientific \&
Technical, 1994,  pp. 46--128}
\lit Pd P1 {P. Podle\'s, \Quote{Quantum Spheres,} Lett. Math. Phys. {\bf 14} (1987),
 193--202}
\lit Pdiii P2 {P. Podle\'s, \Quote{Complex quantum groups and their real
representations,} Publ. RIMS, Kyoto {\bf 28} (1992),  709--745}
\lit Pdii P3 {P. Podle\'s, \Quote{Symmetries of Quantum Spaces. Subgroups and Quotient
 Spaces of Quantum $SU(2)$ and $SO(3)$ Groups,} Commun. 
 Math. Phys., {\bf 170} (1995),  1--20}
\lit PdIV P4 {P. Podle\'s, \Quote{The Dirac operator and gamma matrices 
for quantum Min\-kowski spaces}, q-alg/9703014, will appear in 
J. Math. Phys.}
\lit PW PW {B. Parshall, J. P. Wang, \Quote{Quantum linear groups,} Mem. AMS
{\bf 439}, 1991, AMS}
\lit WP PW1 {P. Podle\'s, S. L. Woronowicz, \Quote{Quantum Deformation of
Lorentz group,} Commun. Math. Phys. {\bf 130} (1990),  381--431}
\lit Pdvi PW2 {P. Podle\'s, S. L. Woronowicz,  \Quote{On the classification of Quantum Poincar\'e Groups,} Commun.
Math. Phys. {\bf 178} (1996),  61--82 and references therein}
\lit Ross R {M. Rosso, \Quote{Alg\`ebres enveloppantes quantifi\'ees, groupes quantiques compacts
de matrices et calcul diff\'erentiel non commutatif,} Duke Math.
J. {\bf 61},1 (1990),  11--40}
\lit FRT RTF {N. Yu. Reshetikhin, L. A. Takhtadzyan%
%
%\foot{name sometimes transcripted to  ``Takhtajan''}
, L. D. Faddeev,
\Quote{Quantization of Lie groups and Lie algebras,} Leningrad Math.
J. {\bf 1}(1) (1990),  193--225}
\lit Schm S {J. Apel, K. Schm\"udgen, \Quote{Classification of
three-dimensional covariant differential calculi on Podles'
quantum spheres and on related spaces,} Lett. Math. Phys.
{\bf 32} (1994),  25--36}
\lit Ta T1 {M. Takeuchi, \Quote{Quantum orthogonal and symplectic groups and their
 embedding into quantum GL,} Proc. Japan Acad., 
 {\bf 65} (1989), Series A, No.~2,  55--58}
\lit Tak T2 {M. Takeuchi, \Quote{Some Topics on $GL_q(n)$,} J. Algebra
{\bf  147} (1992),  379--410}
\lit Wiii W1 {S. L. Woronowicz, \Quote{Twisted $SU(2)$ group. An example of a non-commutative
differential calculus,} Publ. RIMS, Kyoto {\bf23} (1987),  117--181}
\lit Wi W2 {S. L. Woronowicz, \Quote{Compact matrix pseudogroups,} Commun.
Math. Phys. {\bf 111} (1987),  613--665}
\lit TK W3 {S. L. Woronowicz, \Quote{Tannaka-Krein duality for compact matrix
quantum groups. Twisted $SU(N)$ groups,} Invent. Math. {\bf93} (1988),  35--76}
\lit Wii W4 {S. L. Woronowicz, \Quote{New quantum deformation of $SL(2, {\Bbb C})$.
Hopf algebra level,} Rep.  Math. Phys. {\bf30} (1991),  259--269;
cf. also Yu. I. Manin, Le\c cons Coll\'ege de France (1989);  E. E. Demidov, 
Yu.~I.~Manin, E. E. Mukhin, D. V. Zhdanovich,
\Quote{Non-standard quantum deformations of $GL(n)$ and constant solutions
of the Yang-Baxter equation,} Prog. Theor. Phys. Suppl. {\bf102} (1990),
 203--218; S. Zakrzewski, \Quote{A Hopf star-algebra of polynomials
on the quantum $SL(2,\R)$ for a `unitary' R-matrix,}
Lett. Math. Phys. {\bf 22} (1991),  287--289;
H. Ewen, O. Ogievetsky, J. Wess, \Quote{Quantum matrices
in two dimensions,} Lett. Math. Phys. {\bf 22} (1991),
 297--305}
\lit Wv W5 {S. L. Woronowicz, \Quote{Unbounded elements affiliated with $C^*$-algebras
and non-compact quantum groups,} Commun.  Math. Phys. 
{\bf136} (1991),  399--432; S. L. Woronowicz, \Quote{$C^*$-algebras
generated by unbounded elements,} Rev.  Math.
Phys., Vol. 7, No. 3 (1995),  481--521}
\lit WZi WZ1 {S. L. Woronowicz, S. Zakrzewski, \Quote{Quantum Lorentz group
having Gauss decomposition property,} Publ. RIMS, Kyoto, {\bf28} (1992),  809--824}
\lit WZ WZ2 {S. L. Woronowicz, S. Zakrzewski, \Quote{Quantum deformations of the
Lorentz group. The Hopf $*$-algebra level,} Comp. Math. {\bf90} (1994),  211--243}
\lit Zal Z {S. Zakrzewski, \Quote{Realifications of complex quantum groups,} in
\Quote{Groups and related topics,} R. Gielerak et al. (eds.),
Kluwer Academic Publishers, 1992, pp. 83--100}

\def\FG{{\sl Fun}(G)}
\def\Mor{{\rm Mor}}
\def\me{^{-1}}
\def\um#1#2{u\me_{\ #1#2}}
\def\Pol{{\rm Pol}}
\def\Poly{{\rm Poly}}
\def\PolG{{\rm Pol}(G)}
\def\PolyG{{\rm Poly}(G)}
\def\Bew{\removelastskip\smallbreak \rm Proof: }
% equation instead of Gleichung:
\def\Eqlabel#1#2{\iftestmodus(#1)\else%
\global\advance\Cglei by 1
\Label{#1}\Csection\Csubsection\Cglei{Equation} % Verzeichniseintrag
#2 (\number\Cglei) %Gleichungsnummer anschreiben
\iflabelanzeige\rlap{[#1]}\fi\fi}
\def\Rema#1{\Einl{#1}{Remark}\rm}
\def\Thm#1{\Einl{#1}{Theorem}\rm\sl}
\def\Cor#1{\Einl{#1}{Corollary}\rm\sl}

\def\Im{{\rm Im}}
\def\Ke{{\rm Ker}} % Neudefinition
\def\Esy{E^{\rm sym}}
\def\Easy{E^{\rm asym}}
\def\Eins{{\bf 1}}
\def\Smatrix#1{{\def\cr{aus\noexpand\ueber}\def\atop{aa}
\xdef\mmm{(\noexpand\ueber#1aus)}}\def\ueber##1aa##2aus{{##1\atop ##2}}\mmm}
\def\ua#1#2{{\def\qq{#1}\if\qq.u^\alpha\else u^\alpha_{#1#2}\fi}}
\def\ub#1#2{{\def\qq{#1}\if\qq.u^\beta\else u^\beta_{#1#2}\fi}}
\def\cor{corepresentation\ANH}
\def\cors{corepresentations\ANH}
\def\Ho{{\cal H}_0}
\def\tensorop{\mathrel{\tensor^{\rm op}}}\def\normc#1{\Vert#1\Vert_{C^* }}
\def\normi#1{\Vert#1\Vert_{( .\mid. )}}
\def\ttensor{\mathrel{\widehat\tensor}}
\def\cmqg{compact matrix quantum group\ANH}
\def\cmqgs{compact matrix quantum groups\ANH}
\def \uab#1{{\def\g{#1}\ua..^{\if\g.\else#1#1\fi}\tensor \if\g.\ub..\else\ua..\fi^c}}
\def\em{e_{-1}}

% Fussnoten:
\newcount\fuss
\fuss=0
\long\def\foot#1{\advance\fuss by 1\relax\footnote{$^{\number\fuss}$}{{\rm#1}}}

% BEGINNING OF TEXT
%Einleitungsseite
\vglue 0pt plus 1 fil
\seitenzahlfalse
\centerline{\rmzwei Introduction to Quantum Groups}
\bigskip
\centerline{ 
P. Podle\'s{$\dag$}\footnote{$^\ast$}{Supported by Graduiertenkolleg \Quote{Mathematik im Bereich ihrer Wechselwirkung mit
der Physik,}
Dept. of Mathematics, Munich University and by Polish KBN grant No. 
2 P301 020 07}
 and E. M\"uller{$\ddag$}}
\bigskip
{\sevenrm 
$$\vcenter{\baselineskip=0pt\hsize=0.8\hsize\noindent
$\scriptstyle\dag$ Department of Mathematical Methods in Physics, 
Faculty of Physics, University of Warsaw, Ho\.za 74, 00--682 Warszawa, Poland.
E-mail: podles\char64 fuw.edu.pl

\noindent
$\scriptstyle\ddag$ {Graduiertenkolleg
\Quote{Mathematik im Bereich ihrer Wechselwirkung mit der Physik,}
Department of Mathematics, Munich University, Theresienstra\ss e 39, 80333 M\"unchen, Germany. E-mail: emueller\char64 rz.mathematik.uni-muenchen.de}

}$$}

\bigskip
\Abschnitt{Abstract}
We give an elementary introduction
to the theory of algebraic and topological quantum groups (in the spirit
of S. L. Woronowicz). In particular, we recall the basic facts
from Hopf \hbox{($*$-)} algebra theory, theory of compact (matrix) quantum groups
and the theory of their actions on compact quantum spaces.
We also provide the most important examples, including the
classification of quantum $SL(2)$-groups, their real forms and
quantum spheres. We also consider quantum $SL_q(N)$-groups and
quantum Lorentz groups.

\Abschnitt{Contents}
1. Introduction and physical motivations\dotfill\Seitenref{P1}

2. Polynomials on classical groups of matrices\dotfill\Seitenref{P2}

3. Examples of quantum groups\dotfill\Seitenref{P3}

4. $*$-Structures\dotfill\Seitenref{P5}

5. Compact Hopf $*$-algebras\dotfill\Seitenref{P6}

6. Actions on Quantum Spaces\dotfill\Seitenref{P8}

7. Quantum Lorentz groups\dotfill\Seitenref{P9}

8. References\dotfill\Seitenref{P10}

\vglue 0pt plus 1 fil
\break

\subsection{P1}{Introduction and physical motivations}
\Abschnitt{What are \qgs?}
Let~$G$ be a group in the usual sense, i.~e.~a set satisfying the
group axioms, and~$k$ be a field. With this group one can associate a commutative,
associative~$k$-algebra of functions from~$G$ to~$k$
with pointwise algebra structure, i.~e.~for any two elements~$f$ and~$f'$, for
any scalar $\alpha\in k$, and $g\in G$ we have
$$(f+f')(g):=f(g)+f'(g), \ (\alpha f)(g):=\alpha f(g), \
(ff')(g):=f(g)f'(g)$$
If~$G$ is a topological group, usually only continuous
functions are considered, and for an algebraic group the functions are
normally polynomial functions. These algebras are called ``algebras of functions
on~$G$.'' These algebras inherit some extra structures and axioms for those
structures
from the group structure and its axioms on~$G$. Locally compact groups can
be reconstructed from this algebra.

Now  the algebra is {\sl deformed} or {\sl quantized}, i.~e.~the
 algebra structure is changed
so that the algebra is not commutative any more, but the extra structures
and axioms for them remain the same. This algebra is called \Quote{algebra of
functions on a quantum group}, where  \Quote{quantum group} is just an abstract object \Quote{described}
by the deformed algebra. This process can be summarized as follows:
$$\matrix{\vbox{\hbox{classical group $G$}\hbox{axioms of a group}}
&&\vbox{\hbox{quantum group}\hbox{(abstract object)}}\cr
\downarrow&&\updownarrow\cr
\vbox{\hbox{commutative algebra of}\hbox{functions on $G$ with}
\hbox{corresponding extra axioms}}&
\buildrel\textstyle{\rm forget\  about\atop \rm group}\over\llongrightarrow&
\vbox{\hbox{non-commutative algebra with}
\hbox{same extra axioms; ``algebra of}
\hbox{functions on a quantum group''}}\cr}$$
There is a similar concept of \Quote{quantum spaces}: If~$G$ acts on
a set~$X$ (e.~g.~a vector space), there is a corresponding so-called
{\sl coaction} of the commutative algebra of functions on~$G$ on
the commutative algebra of functions on~$X$ satisfying certain axioms.
The latter algebra can often be deformed/quantized
into a non-commutative algebra,
called the \Quote{algebra of functions on a quantum space} with a
similar coaction.
There are three ways of considering algebras of functions on a group and
their deformations:
\item (a) polynomial functions $\PolyG$ (developed by Woronowicz and
Drinfel'd),
\item (b) continuous functions $C(G)$, if~$G$ is a  topological group
(developed by Woronowicz),
\item (c) formal power series (developed by Drinfel'd).

Only the first two approaches will be dealt with in the sequel. They include
representation theory, Peter-Weyl theory, Tannaka-Krein theory,
and actions on quantum spaces.

There is a second approach to quantum groups. If~$G$ is a connected,
simply connected Lie group, $G$
can be reconstructed from the universal enveloping algebra~$U(\frak g)$
of the corresponding Lie algebra~$\frak g$. The algebra~$U(\frak g)$ again inherits
some extra structures and axioms and can be deformed. The deformed universal
enveloping algebra can be regarded as universal enveloping algebra corresponding
to a quantum group.
One can consider
\item (d) the quantized universal enveloping algebra $U_q(\frak g)$ (developed
 by Jimbo),
\item (e) formal power series (to be more precise, the ring of formal
power series in~$h$ over a free algebra, 
subject to certain relations which are the
same as for~$U(\frak g)$ in the case~$h=0$. From this ring the algebra~$U_q(\frak g)$
can be extracted. This approach has been developed by Drinfel'd).

This approach will not be used in the sequel.

\Abschnitt{Physical motivations}
There are some physical motivations for \qgs including
\item 1. integrable models---handled with approach (e),
\item 2. conformal field theory---handled with approach (e),
\item 3. physical models based on quantized space-time---handled with
approaches (a), (b), and (e).

The last motivation shall be explained in more detail. One of the main
problems in Quantum Field Theory (QFT) is to join QFT and General Relativity
Theory in a consistent way. It seems that in such a new theory it would
be impossible to study the
geometry of the space when very small volumes are considered.
If you
consider a cube in space, each vertex of it having Planck's length or less, and
measure simultaneously the three coordinates $x$, $y$,
and $z$ of a particle in it, then the uncertainty of the measurement, i.~e.~the
errors $\Delta x$, $\Delta y$, and $\Delta z$ are very small, whence by
Heisenberg's uncertainty relation the errors of the coordinates of the momentum
are big and therefore the uncertainty of the energy~$\Delta E$ is big, too.
Since the energy is positive, the expected value~${\langle E\rangle}$
of the energy is big, and the smaller
the cube the bigger the energy, which at a certain stage generates
a black hole.
Therefore the observation of the geometry of the
space gives it a different geometry, which makes this observation useless
(We have used here the arguments by Professor W. Nahm).

Quantum mechanics says that physical quantities such as momentum and position,
which can be measured, correspond to self-adjoint operators on a Hilbert
space. Its elements describe possible states of a physical system. When a
quantity is measured, the state is projected onto an eigenvector of the
operator, and the result of the measurement is the corresponding
eigenvalue. Two quantities can be measured simultaneously if and only if
the corresponding operators commute.
In usual quantum mechanics the operators corresponding to
the three coordinates of space commute and can be measured simultaneously,
which leads to the problem with the black hole.
Thus it is reasonable to assume that the operators corresponding to
the coordinates $x$, $y$, and $z$
do not commute (whence they cannot be measured simultaneously).
Hence the commutative algebra generated by the operators corresponding to
$x$, $y$, and $z$, which is isomorphic to the algebra of
polynomials on $\R^3$, is replaced by a non-commutative algebra on a
quantum space. In order to give sense to self-adjoint operators,
this algebra should be a $*$-algebra.
\Def{P1.1} \voritem (a) A {\sl $*$-algebra} is a \C-algebra~$A$ equipped with an antilinear,
antimultiplicative involution $*\colon A\to A$, i.~e.~for all $a,b\in A$
and $\lambda\in \C$ the following holds:
$$(a+b)^*=a^*+b^*, \ (\lambda a)^*=\bar\lambda a^*, \
 (ab)^*=b^*a^*, \ (a^* )^*= a.$$
\item (b) Let $A$, $B$ be $*$-algebras. An algebra homomorphism $\phi\colon A\to B$
is called {\sl $*$-homo\-morphism}, if $\phi(a^* ) = \phi(a)^*$ for all $a\in A$.
\Ende
Physical experiments should be comparable and reproducible, i.~e.~the
same experiment performed at different places and times ought to give the
same result. Therefore the theory should be invariant with respect to
certain symmetry groups (containing translations in time and
space). But the classical (symmetry)
groups do not
fit well to quantum spaces, so they have to be changed to quantum groups,
too. (Example: The group $SO_3(\R)$ of rotations in three-dimensional
space acts on the sphere~$S^2$. When the algebra of functions on~$S^2$
is properly deformed such that the algebra
becomes non-commutative, then there is
no reasonable coaction of the usual algebra of functions on $SO_3(\R)$
any more. \Pd{Remark~2})

There is another motivation---deformation of an existing physical theory 
may help to
understand the theory in a better way. It can reveal why the theory works, 
what is a consequence, and what
is just a coincidence.

Example \PdIV{}: 
After looking at deformations of standard Dirac theory,
the covariance
of the Dirac equation can be seen more directly---on the level of groups
rather than Lie algebras. 
For the wave vector~$\Psi$ there is the
equation $\bar\Psi=\Psi^{\scriptstyle\dag} \gamma^0$, where $\gamma^0$ also appears
in the Dirac equation. In the deformed theory there is $\bar\Psi=\Psi^{\dag} A$
with $A\ne \gamma^0$ in general, 
so that $A=\gamma^0$ is just a coincidence,
and the condition $A=\gamma^0$ is not really important for the
theory.

In physics all symmetry groups are groups of matrices or can be described with
groups of matrices, therefore the case of matrix groups is considered.

\Abschnitt{Acknowledgment}
These lecture notes were written down by E. M. after the
lectures by P. P.  given at the Department of Mathematics,
Munich University. The first author is very grateful to Professor Dr. Hans-J\"urgen Schneider for his warm hospitality at 
Munich University. We thank him for his useful remarks.
An earlier version of these lectures was given by P. P.
at Kyoto University in 1990--91. The first author would
like to express his gratitude to Professor Huzihiro Araki for his
kind hospitality and encouragement.

\subsection{P2}{Polynomials on classical groups of matrices}

\Abschnitt{Notations}

In the sequel the base field of all vector spaces and algebras
is the field \C of complex numbers. A {\sl unital algebra} is an
(associative) algebra with a unit element, and a {\sl unital mapping} is a mapping
between unital algebras
which sends the unit element to the unit element.

Let \N, $\N_0$ and \R denote the sets of positive integers, non-negative
integers and 
real numbers respectively and fix $M,N\in\N$.
Let $A$ be a unital algebra and let $M_{M\times N}(A)$ be the vector space of
$M\times N$-matrices with entries in~$A$. If $M=N$, $M_N(A):=M_{N\times N}(A)$
is a unital algebra. For each matrix
$M\in M_{M\times N}(A)$ let $M_{ij}$ be
the entry at the $i$-th row and $j$-th column of $M$. 
Let~$B$ be another algebra, $\phi\colon A\to B$ a map and $M\in M_{M\times N}(A)$.
Then $\phi(M)$ is shorthand for the matrix in~$M_{M\times N}(B)$
with entries~$\phi(M_{ij})$.
The group $GL(N,\C)$ of invertible
$N\times N$-matrices with complex entries is equipped with a topology
inherited from the norm topology of the vector space $M_N(\C)\cong \C^{N^2}$.
The neutral element of a group is denoted by~$e$.

Let $\C^N$ denote the space of row vectors and~${}^N\C$ the space
of column vectors. Using matrix multiplication, $\C^N$ can be regarded as dual space
of ${}^N\C$. If $\{e_1, \ldots, e_N\}$ is a basis of~${}^N\C$,
then there is a dual basis $\{e_1', \ldots, e_N'\}$ of $\C^N$ such
that $e_i'e_j = \delta_{ij}$ for all $i,j\le N$. In a similar way there
are dual bases of the $k$-fold tensor products~$({}^N\C)^{\tensor k}$
and~$(\C^N)^{\tensor k}$.

In the sequel
the indices $i$, $j$, $i'$, $j'$,
$k$ denote positive integers less or equal to~$N$.

Let $\Eins_N$ denote the identity matrix with~$N$ rows and columns or
the identity endomorphism of~$\C^N$ or~${}^N\C$.

\Abschnitt{Functions on groups}

Let $G$ be an arbitrary subgroup of the group $GL(N,\C)$.
Let $\FG$ be the algebra of complex valued functions on~$G$.
This algebra is unital with unit element
$\underline 1\colon G\to \C, \ g\mapsto 1$ and  is a $*$-algebra, where
for all $f\in \FG$ the function $f^*$ is defined by $f^*(g):=\overline{f(g)}$
for all $g\in G$.

For all~$i$ and~$j$, the
coefficient functions
$$u_{ij}\colon G\to \C, \ g\mapsto g_{ij}\hbox{ and }
  \um ij \colon G\to \C, \ g\mapsto (g\me)_{ij}$$
belong to $\FG$.
Then the matrices $u:= (u_{ij})_{1\le i,j\le N}$ and
$u\me:=(\um ij)_{1\le i,j\le N}$ belong to $M_N(\FG)$ and
are inverses of each other in $M_N(\FG)$. This justifies the notation~$u\me$.
\Def{P2.1} Let $\PolG$ be the subalgebra of $\FG$
generated by the elements $u_{ij}$ and $\um ij$ for all $i$ and~$j$.
\Ende
Remark: This algebra is automatically unital because of the relation
$\underline 1 = \sum_{k=1}^n u_{1k}\um k1$.
The algebra is called \Quote{algebra of
holomorphic polynomials on~$G$}, too.

\Lemma{P2.2} If $G\subset SL(N,\C)$ then $\PolG$ is already
generated by the elements~$u_{ij}$.
\Bew By the usual formula for  the inverse of a matrix,
$(g\me)_{ij} = (-1)^{i+j}\det{\tilde g_{j,i}}/\det(g)$ for
all $g\in G$, where
the $(N-1)\times (N-1)$-matrix $\tilde g_{j,i}$
is obtained from~$g$ by deleting the
$j$-th row and the $i$-th column. But $\det(g)=1$, whence also $\um ij$
is a polynomial in the functions $u_{i'j'}$.
\Ende

\Def{P2.3} Let $\PolyG$ be the $*$-subalgebra of $\FG$
generated by the elements $u_{ij}$ and $\um ij$.
\Ende
Usually the algebra $\PolyG$ is considerably bigger than $\PolG$.

\Lemma{P2.4} If $G$ is a compact subgroup of $GL(N,\C)$, then
$\PolyG$ is generated by  the elements $u_{ij}$ as $*$-subalgebra{}.
\Bew The map $\phi\colon G\to \R^+$, $g\mapsto \vert\det(g)\vert$ is
a group homomorphism from~$G$ into the multiplicative group of positive
real numbers. Since $\phi$ is continuous and~$G$ is compact, the image
of~$\phi$ is a compact subgroup of~$\R^+$. But $\{1\}$ is
the only compact subgroup of~$\R^+$, whence $\phi(g)=1$ for all
$g\in G$. Therefore
$$1=\det(g)\overline{\det(g)}= \det(g)\det((\bar g_{ij})_{1\le i,j\le N}).$$
Thus $\det(u)$ is invertible in $\PolyG$ with inverse
$\det( (u^*_{ij})_{1\le i,j\le N} )$, whence the elements $\um ij$
can be expressed by the $u_{i'j'}$ and $u^*_{i'j'}$.
\Ende

\Rema{P2.5} Let $I$ be an index set and let $G$ be a subgroup of
$\prod_{\alpha\in I} GL(N_\alpha,\C)$. Each element~$g$ of
this group can be written as $g=(g_\alpha)_{\alpha\in I}$ with
$g_\alpha\in GL(N_\alpha,\C)$ for all $\alpha\in I$ and define
$$u^\alpha_{ij}, (u^{\alpha})^{-1}{}_{ij}
 \colon G\to \C, \
u^\alpha_{ij}(g):=(g_\alpha)_{ij}, \
  (u^{\alpha})^{-1}{}_{ij}:=(g\me_\alpha)_{ij}$$
for all $g\in G$. The algebras $\PolG$ and $\PolyG$ are generated
by the elements $u^\alpha_{ij}$ and $(u^{\alpha})^{-1}{}_{ij}$
as algebras or $*$-algebras, respectively.
This generalization covers all compact groups~$G$, because the
group homomorphism
$$G\to \prod_{\pi\in \widehat G} GL(\dim(\pi),\C), \
g\mapsto (\pi(g))_{\pi\in\widehat G}$$
where $\widehat G$ is the set of finite dimensional irreducible representations
of~$G$, is injective if $G$ is compact (cf. Tannaka-Krein duality).
\Ende

The multiplication, unit, and the inverse on~$G$ lead to the following extra
structures on $\FG$:

$$\eqalign{\Delta\colon \FG\to {\sl Fun}(G\times G),& \quad
(\Delta f) (g,h):=f(gh) \hbox{ for all $g,h\in G$ (Comultiplication)}, \cr
 \eps\colon \FG\to \C,& \quad
 \eps (f):=f(e) \hbox{ (Counit),} \cr
 S\colon \FG\to \FG,& \quad (Sf)(g):=f(g\me) \hbox{ for all $g\in G$
 (Antipode).} \cr}$$
These maps are unital $*$-homomorphisms.
The (algebraic) tensor product $\FG\tensor\FG$ is 
the vector subspace of ${\sl Fun}(G\times G)$ generated by elements
$u\tensor v$, where $u,v\in \FG$, by defining
$(u \tensor v)(g,h):=u(g)v(h)$ for all $g,h\in G$. Equality only holds if~$G$
is finite.

The axioms for the group structure on~$G$ are reflected by certain axioms
for the extra structures on $\FG$. Let~$f$ be an element of~$\FG$
such that $\Delta(f)\in \FG\tensor\FG$.
Since the multiplication in~$G$ is associative, we have
$$(\Delta\tensor\id)\Delta ( f)=(\id\tensor\Delta)\Delta (f). \Eqno 1$$
The property of the neutral element, namely $ge=eg=g$ for all $g\in G$, leads to the
equation
$$(\eps\tensor\id)\Delta(f)=(\id\tensor \eps)\Delta(f)=f.\Eqno 2 $$
(Here the usual identification $\C\tensor V\cong V\tensor \C\cong V$ for
all \C-vector spaces is used).
Let the linear map $\mu\colon \FG\tensor\FG\to \FG$, $f\tensor f'\to ff'$ be induced
by the multiplication in $\FG$. Then the properties $gg\me=g\me g=e$ of the inverse can
be expressed as
$$\mu(S\tensor \id)\Delta(f) = \mu(\id\tensor S)\Delta(f)
= \eps(f)\underline 1.\Eqno3$$
\Def{P2.6} A unital algebra~$H$ is called {\sl Hopf algebra}, if
there are unital algebra homomorphisms $\Delta\colon H\to H\tensor H$
and $\eps\colon H\to \C$
and a linear map $S\colon H\to H$
satisfying axioms (\SRefg{1})--(\SRefg{3}) for all $f\in H$.
\Ende

The following lemma gives examples of Hopf algebras and shows why $\PolG$
and the elements $\um ij$ are interesting.
\Lemma{P2.6a} $\PolG$ is a Hopf algebra satisfying
$$\displaylines{\Delta u_{ij} = \sum_{k=1}^N u_{ik}\tensor u_{kj},  \
\Delta \um ij = \sum_{k=1}^N \um kj\tensor\um ik, \cr
\eps(u_{ij})=\eps(\um ij)=\delta_{i,j}, \
S(u_{ij})= \um ij, \ S(\um ij)= u_{ij}.}$$
If~$G$ is finite, also $\FG$ is a Hopf algebra.
\Bew For all $g,h\in G$,
$$\Delta u_{ij} (g,h) = u_{ij}(gh) = (gh)_{ij} =
\sum_{k=1}^N g_{ik}h_{kj} = \sum_{k=1}^N u_{ik}(g)u_{kj}(h) =
\sum_{k=1}^N (u_{ik}\tensor u_{kj})(g,h).$$
A similar computation yields the formula for $\Delta(\um ij)$.
Therefore the image of~$\PolG$ under $\Delta$ is contained in $\PolG\tensor \PolG$.
The values of the counit can be computed:
$\eps(u_{ij})=\eps(\um ij) = e_{ij} = \delta_{i,j}$.
The equations for the antipode follow from
$(S(u_{ij}))(g) = u_{ij}(g\me) = \um ij (g)$ for all $g\in G$.
The Hopf algebra
axioms are clearly satisfied, because $\PolG$ is a subalgebra of $\FG$.

If~$G$ is finite, then
$\FG$ is a Hopf algebra because $\FG\tensor\FG= {\sl Fun}(G\times G)$.
\Ende

The following general theorem for Hopf algebras can be inferred from \Abe{}.
\Thm{P2.7} Let $H$ be a Hopf algebra with unit element~1.
\item (a) The maps $\eps$ and $S$ are unique if $\Delta$ is fixed.
\item (b) $S$ is a unital antihomomorphism
\item (c) If $\tau\colon H\tensor H\to H\tensor H$, $x\tensor y\mapsto
y\tensor x$ denotes the flip automorphism, then
$$\Delta S = \tau (S\tensor S)\Delta, \ \eps S=\eps.$$
\item (d) Let  $S'\colon H\to H$ be a \C-linear map. Then the following
are equivalent:
\itemitem (i) $\mu(\id\tensor S')\tau\Delta(f) = \mu(S'\tensor \id)\tau\Delta(f)=
\eps(f)1$ for all $f\in H$,
\itemitem (ii) $S\circ S'=S'\circ S=\id$.
\Ende

\Rema{P2.8} \voritem (a) In general, the antipode of a Hopf algebra is not
invertible.
\item (b) A map $S'$ such as in part (d) of \Ref{P2.7} is called
{\sl skew antipode}, and there is another Hopf algebra structure on~$H$
with comultiplication $\tau\Delta$, counit~$\eps$ and antipode~$S'$.
\item (c) A motivation for the fact, that the counit, but not the antipode is
 an algebra homomorphism, if~$H$ is not commutative:
 Since $\Delta$ and the identity are algebra homomorphisms, there is no reason
 following from axiom (\SRefg{2}) that $\eps$ should not be an algebra homomorphism.
 But the map $\mu$ in axiom (\SRefg{3}) is an algebra homomorphism if and only if
 $H$ is commutative. Therefore it should not be expected that~$S$ is an
 algebra homomorphism.
\Ende

For all $f\in \FG$ satisfying $\Delta(f)\in \FG\tensor \FG$,
the following equation holds:
$$\eqalign{\Delta(f^* )(x,y) &= f^* (xy) = \overline{f(xy)} =
\overline{\Delta f(x,y)} = \overline{\sum f_1(x)f_2(y)} =\cr
&=
\sum f_1^*(x)f_2^*(y) = ( *\tensor * )\Delta(f) (x,y)}$$
for all $x,y\in G$. This motivates the following definition.

\Def{P2.9} A unital algebra $H$ is called a {\sl Hopf $*$-algebra}, if $H$ is
both a Hopf algebra and a $*$-algebra such that
$\Delta(f^* ) = ( *\tensor * )\Delta f$ for all $f\in H$.
\Ende

{}From the definitions and \Ref{P2.6a} follows immediately

\Lemma{P2.9a} $\PolyG$ is a Hopf $*$-algebra, and if $G$ is finite, also
$\FG$ is a Hopf $*$-algebra.
\Ende

\Prop{P2.10} Let $H$ be a Hopf $*$-algebra. Then
\item (a) For all $x\in H$, $\eps(x^* ) = \overline{\eps(x)}$, i.~e.~$\eps$
is a $*$-homomorphism.
\item (b) $S\circ *\circ S\circ * = \id$, in particular, $S$ is bijective.
\Bew \voritem (a) Since the map $H\to \C$, $x\mapsto \overline{\eps(x^* )}$
satisfies the properties of the counit, both are equal by \Ref{P2.7}, part~(a),
whence the assertion follows.
\item (b) The map $*\circ S\circ *$ satisfies all properties of the
skew antipode. By \Ref{P2.7}, part (d) it is equal to it. This implies the
two equivalent equalities
$*\circ S\circ *\circ S = \id_H = S\circ *\circ S\circ *$.
\Ende

\Abschnitt{Elements of representation theory}
Let $H$ be a Hopf algebra.
\Def{P2.11} Let $k$ be a positive integer.
A matrix $v\in M_k(H)$ is called {\sl \cor\empty}, if the entries
satisfy the following relations for all indices $a$ and $b$.
\item (a) $\Delta v_{ab} = \sum_{c=1}^k v_{ac}\tensor v_{cb}$,
\item (b) $\eps(v_{ab})=\delta_{a,b}$

The number $\dim v := k$ is called the {\sl degree} of the \cor,
and the elements $v_{ab}$ are called the {\sl matrix elements} of the
\cor.
\Ende

\Rema{P2.12} \voritem (a) Let $v$ be a \cor of a Hopf algebra~$H$.
Then $S(v_{ab})=(v\me)_{ab}$ for all indices $a$, $b$. Thus Condition~(b)
of \Ref{P2.11} can be equivalently replaced by invertibility of~$v$
(note that Condition~(a) implies ${\eps(v)v=v}$).
\item (b) Let~$G$ be a classical group of matrices and $H$ one of
the Hopf algebras $\PolG$ or $\PolyG$. Let
$$v\colon G\to M_k(\C), \ g\mapsto (v_{ab}(g))_{1\le a, b\le k}$$
be a map such that all functions $v_{ab}$ are contained in~$H$. Then
$(v_{ab})_{1\le a,b\le k}$ is a \cor if and only if
$v$ is a representation of~$G$.
\Bew \voritem (a) This follows from the axioms for the antipode of a
Hopf algebra.
\item (b) For all $x,y\in G$ the following equations hold.
$$\displaylines{(\Delta v_{ab})(x,y)=v_{ab}(xy)=(v(xy))_{ab}, \cr
(\sum_{c=1}^k v_{ac}\tensor v_{cb})(x,y) =
\sum_{c=1}^k v_{ac}(x)v_{cb}(y) = (v(x)v(y))_{ab}.}$$
Therefore condition (a) in \Ref{P2.11} is equivalent to $v(xy)=v(x)v(y)$.
A computation of $\eps(v_{ab})$ shows that condition (b) is equivalent
to $v(e)=\Eins_k$.
\Ende
\goodbreak
Now fix a Hopf algebra~$H$.
\Def{P2.13}
Let $v$ and $w$ be two \cors of $H$.
\item (a) Then
$v\oplus w$ and $v\otimes w$ are \cors of~$H$, where $v\oplus w$
is a matrix with $\dim(v)+\dim(w)$ rows and columns given by
$$\pmatrix{v&0\cr 0&w},$$
and the matrix of $v\otimes w$ has $\dim(v)\dim(w)$ rows and columns and
entries given by $(v\otimes w)_{ij,kl}:= v_{ik}w_{jl}$, where
the indices $i$, $k$ take values between 1 and $\dim(v)$ and the
indices $j$, $l$ between 1 and $\dim(w)$.
\item (b) A  $\dim(w)\times\dim(v)$ matrix
$A$ over \C {\sl intertwines} $v$ with $w$, if $Av=wA$.
Define $\Mor(v,w)$ as vector space of intertwining matrices between $v$
and~$w$. The elements of $\Mor(v,w)$ can be regarded as \C-linear
maps from $\C^{\dim v}$ to $\C^{\dim w}$. The \cors $v$ and $w$ are said
to be {\sl equivalent} ($v\cong w$) if $\dim(v)=\dim(w)$ and there is an invertible
element in $\Mor(v,w)$.
\Ende
\DefLemma{P2.14m} Let $w$ be a \cor of dimension~$N$
and $V\subset {}^N\C$ a subspace of dimension~$l$.
Then the following are equivalent:
\item (a) For each $\rho\in \Hom(H, \C)$ the statement $\rho(w)V\subseteq V$ holds.
\item (b) There is a \cor $v$
and a basis $a_1, \ldots, a_l$ of $V$ such that for the $N\times l$-matrix
${A_l:=(a_1\cdots a_l)}$
the equation $wA_l = A_lv$ holds. This is equivalent to
the condition that~$A_l$ is an injective intertwiner of~$v$ with~$w$.
\item (c) There is a \cor $v$ of dimension $l$ and an invertible matrix~$A$,
the first $l$ columns of which are a basis of~$V$ and such that
$$wA = A\pmatrix{v&*\cr 0&* }.$$

If one of the equivalent conditions holds, then $V$ is called \Quote{$w$-invariant
subspace}, and the \cor $v$ in part (b) and (c) is called \Quote{sub\cor
of $w$} and we write $v=w\vert_V$ (Note that $v$ depends on the chosen basis
of~$V$).

\Bew \voritem (a)$\Folgt$((b)$\iff$(c)).  Let $a_1, \ldots, a_l$ be a basis of~$V$ and extend
it to a basis $a_1, \ldots, a_N$ of ${}^N\C$. Then let $A_l$ be the
$N\times l$-matrix $(a_1\cdots a_l)$ and~$A$ be the $N\times N$-matrix
$(a_1\cdots a_N)$.
Then~$A$ is invertible and let $B:=A\me w A$.
Let $\rho\in \Hom(H, \C)$.
Then $A\me \rho(w) A= \rho(B)$.
Now condition~(a) means that there is a matrix $C_\rho\in M_l(\C)$
such that $\rho(w)A_l=A_l C_\rho$, whence $\rho(B)$ looks like
$$\pmatrix{\rho(v)&*\cr 0&* },$$
where $v$ is the submatrix of~$B$ consisting of the first~$l$ rows and
columns.
Since this holds for all linear forms, there is the matrix equation
$$A\me w A= \pmatrix{v& *\cr 0&* } \iff wA=A\pmatrix{v& *\cr 0&* }$$
or, equivalently, by restriction $wA_l=A_l v$.
\item (b) $\!\!\Folgt$(a). From (b) it follows for all $\rho\in\Hom(H, \C)$
that $\rho(w)A_l = \rho(v)A_l$, which gives ${\rho(w)V\subseteq V}$.
\Ende

\Def{P2.14z} Let $w$ a \cor.
\item (a) $w$ is said to be {\sl irreducible} if $w\ne 0$ and there is no sub\cor $v$ such that $0<\dim(v)<\dim(w)$.
\item (b) $w$ is called {\sl completely reducible} if $w$ is equivalent
to a direct sum of irreducible sub\cors.
\Ende

\Lemma{P2.14n} The intersection of invariant subspaces is an invariant subspace.
\Bew This follows directly from  \Ref{P2.14m}, part (a).
\Ende

\Lemma{P2.14o}  Let $A\in \Mor(v,w)$. Then $\Ke(A)$ is $v$-invariant
and $\Im(A)$ is $w$-invariant.
\Bew Use \Ref{P2.14m}, part (a).
For each $\rho\in\Hom(H,\C)$ the equation $A\rho(v)=\rho(w)A$
follows. If $x\in \Ke(A)$ then $A\rho(v)x=\rho(w)Ax=0$, whence $\rho(v)x\in\Ke(A)$ and the kernel is $v$-invariant.
If $y\in\Im(A)$, say $y=Az$, then $\rho(w)y=\rho(w)Az=A\rho(v)z$ is in the
image of~$A$, too.
\Ende
\doppelpunktfalse\Lemma{P2.14} (Schur): Let $v$, $w$ be irreducible
\cors. If $v$ and $w$ are not equivalent, then $\Mor(v,w)=\{0\}$.
If $v$ is irreducible, then $\Mor(v,v)=\C \hbox{{\bf 1}}$,
where {\bf 1} is the identity.
\Bew Let $A\in \Mor(v,w)\setminus\{0\}$. Since $v$ and~$w$
are irreducible, by \Ref{P2.14o}, $A$ must
be injective and surjective, whence $v$ and~$w$ are equivalent. Now
let $w=v$ and~$\lambda$ be an eigenvalue of~$A\in\Mor(v,v)$. Then
${A-\lambda\Eins\in \Mor(v,v)}$ is not
injective and therefore vanishes.
\Ende

\Rema{P2.14a} 
There is a relationship between finite dimensional
right comodules of~$H$ and \cors.
\Ende

\Thm{P2.15} Let $H$ be a Hopf algebra.
\item (a) The matrix elements of \cors
span~$H$.\foot{This result is related to the fact that each
element of a Hopf algebra is contained in a finite dimensional
subcoalgebra.}
\item (b) The matrix elements of a set of non-equivalent irreducible \cors are linearly
independent.
\item (c) The following are equivalent:
\itemitem (i) There is a set $T$ of non-equivalent irreducible \cors
such that the matrix elements of them form a basis of $H$.
\itemitem (ii) Each \cor is completely reducible.\foot{In the language
of Hopf algebras this means that~$H$ is cosemisimple.}
\item {} Moreover if (i) holds then $T$
contains all non-equivalent irreducible \cors.
\Bew \voritem (a) Let $x\in H$. Then there is a number $N\in \N$,
linearly independent elements $x_1, \ldots, x_N$ and $y_1, \ldots, y_N$ in~$H$ such that
${\Delta(x) = \sum_{j=1}^N x_j\tensor y_j}$. By coassociativity,
$\sum_{j=1}^N \Delta(x_j)\tensor y_j = \sum_{j=1}^N x_j\tensor \Delta(y_j)$,
whence there are elements $v_{ij}$ of $H$ such that
$$\Delta(x_j) = \sum_{i=1}^N x_i\tensor v_{ij}$$
for all $j$. Using coassociativity and the properties
for the counit, from these equations it follows that the elements~$v_{ij}$
are matrix elements of a
\cor and $x_j=\sum_{i} \eps(x_i)v_{ij}$ for all~$j$. But then
$$x = \sum_{j=1}^N x_j\eps(y_j) = \sum_{i,j=1}^N \eps(y_j)\eps(x_i)v_{ij}$$
is a linear combination of matrix elements.
\item (b) Use the arguments in the proof of \Wi{Proposition 4.7}.
\item (c)   The conclusion (ii)$\Folgt$(i) is now obvious, because by (a),
the Hopf algebra is spanned by matrix elements of irreducible \cors,
which are linearly independent by (b).
The conclusion (i)$\Folgt$(ii) is proved in \Pdiii{Appendix}. The last
remark follows from~(b).
\Ende
\Prop{P2.16} Let $\{v_\alpha\mid \alpha\in I\}$ and
$\{v'_\beta\mid\beta\in J\}$
be sets of irreducible \cors such that
$$\bigoplus_{\alpha\in I} v_\alpha \cong \bigoplus_{\beta\in J}
v'_\beta.$$
Then the multiplicities of  equivalence classes of irreducible
\cors are the same on both sides.\foot{cf. Krull-Remak-Schmidt theorem}
\Bew The set $\Mor(\bigoplus v_\alpha, \bigoplus v'_\beta)$
can be computed using Schur's lemma (\Ref{P2.14}). But this set
must contain an invertible element, since both direct sums are equivalent.
\Ende
\DefLemma{P2.17} \voritem (a) Let $w$ be a \cor of a
Hopf algebra~$H$. Then also the matrix
$w^c$ with matrix elements $w^c_{ij} := S(w_{ji})$ is a \cor,
the {\rm contragradient \cor\empty} to~$w$.
\item (b) Let $w$ be a \cor of a Hopf $*$-algebra~$H$.
Then also the matrix $\bar w$ with matrix elements $\bar w_{ij}:=w^*_{ij}$
is a \cor. Define~$w^*$ to be the transpose of~$\bar w$.
\item (c) A \cor~$w$ of a Hopf $*$-algebra is called
{\sl unitary} if $\bar w = w^c$ or equivalently $ww^*=w^*w=\Eins_{\dim w}$.
\Bew (a) and (b) follow from the identities $\Delta\circ S = \tau(S\tensor S)\Delta$,
$\eps S = \eps$, $\Delta\circ * = ( *\tensor* )\Delta$.
\Ende

\subsection{P3}{Examples of \qgs}
\Abschnitt{Quantum $SL(2)$-groups}
The simplest Lie group over the complex numbers, which is
interesting and important in physics, is $SL(2,\C)$.
We want to find quantum analogues of ${\rm Pol}(SL(2,\C))$.
The \cors of this Hopf algebra have the following properties:
\item (1) The irreducible \cors are $w^\alpha$, where
$2\alpha\in\N_0$.
\item (2) $\dim(w^\alpha) = 2\alpha+1$ for all $\alpha$,
\item (3) $w^\alpha\tensor w^\beta \cong w^{\vert \alpha-\beta\vert}\oplus
w^{\vert \alpha-\beta\vert+1}\oplus\cdots w^{\alpha+\beta}$
(Clebsch Gordan),
\item (4) Each \cor is completely reducible, or equivalently,
the matrix elements $w_{ij}^\alpha$ span the Hopf algebra.

Remark: The fundamental \cor is $w:=w^{1/2}$ given by
$$g\mapsto (g_{ij})_{1\le i,j\le 2}$$
for $g\in SL(2,\C)$, and $w^0$ is the identity.

\Def{P3.1} A {\sl quantum $SL(2)$-group} is a Hopf algebra satisfying
the properties \hbox{(\SRefg{1})--(\SRefg{4})}.
\Ende
\Thm{P3.2} Up to isomorphism  there are the following quantum $SL(2)$-groups $\cal H$.
The Hopf algebra $\cal H$ is generated by the matrix elements $w_{ij}$
($1\le i,j\le 2$)
of the fundamental \cor $w:=w^{1/2}$ and relations
$$(w\tensor w)E=E, E'(w\tensor w)=E',$$
where the base field~\C is canonically embedded into~$\cal H$ and
there is the following extra relation between the
row vector $E'\in \C^2\tensor\C^2$ and the column vector
$E\in {}^2\C\tensor{}^2\C$: Let $\{e_1, e_2\}$ be a basis of $^2\C$
and $\{e_1', e_2'\}$ be a dual basis of $\C^2$. There is the following
presentation:
$$E=\sum_{i,j=1}^2 E_{ij}e_i\tensor e_j,
  E'=\sum_{i,j=1}^2 E'_{ij}e'_i\tensor e_j'.$$
Then the $2\times 2$ matrices with entries $E_{ij}$ and $E'_{ij}$
are inverses.
There is a basis $\{e_1, e_2\}$ of $^2\C$ such that
$$E=e_1\tensor e_2-qe_2\tensor e_1 \widehat= \pmatrix{0&1\cr-q&0}
\hbox{ or } E= e_1\tensor e_2-e_2\tensor e_1+e_1\tensor e_1
\widehat= \pmatrix{1&1\cr -1&0},$$
where $q\in\C\setminus\{0\}$ must not be a non-real root of unity.
In the first case the quantum group is called the standard deformation~$SL_q(2)$,
in the second case it is called the non-standard deformation $SL_{t=1}(2)$.
The non-standard deformation $SL_{t=1}(2)$ is not isomorphic to any of
the standard deformations, and two standard deformations~$SL_q(2)$ and
$SL_{q'}(2)$ are isomorphic if and only if $q=q'$ or $qq'=1$.
\Ende
\Rema{P3.3} \voritem (a) There is a set of non-standard deformations $SL_t(2)$
indexed by a parameter $t\in\C\setminus\{0\}$ corresponding to the
vector
$E_t=e_1\tensor e_2-e_2\tensor e_1+t e_1\tensor e_1$, but they are all
equivalent to the deformation for $t=1$, because if the basis vector
$e_1$ is replaced by $e_1'=e_1 t$ then
$$t E_t = e_1'\tensor e_2-e_2\tensor e_1'+e_1'\tensor e_1'\widehat= E_1.$$
Since the relations remain the same when~$E$ is multiplied by a non-zero
scalar, the Hopf algebras are isomorphic.
\item (b) For $t\to 0$, the vector $E_t$ tends to the vector for $q=1$.
\item (c) Parts of the proof of \Ref{P3.2} can be found e.~g.~in
\DV{}, \Wii{}, \KP{}.
%The proof that there are no extra relations in the quantum $SL(2)$-groups
%and the treatment on isomorphism classes is not contained therein.
\Ende

To prepare the proof of \Ref{P3.2}, some extra definitions and lemmas
are useful.
\Def{P3.3a} Let $q$ be a complex number. Then a {\sl Hecke algebra}
of degree~$n$ is a unital algebra generated by elements $\sigma_1, \ldots,
\sigma_{n-1}$ subject to the relations
$$\eqalign{ \sigma_k\sigma_{k+1}\sigma_k &= \sigma_{k+1}\sigma_k
\sigma_{k+1} \hbox{ for $1\le k\le n-2$,} \cr
(\sigma_k-1)(\sigma_k+q^2 1) &=0, \cr
\sigma_k\sigma_l&=\sigma_l\sigma_k \hbox{ for $\vert k-l\vert\ge 2$}.}$$
\Ende

{}From these relations follows an important property of Hecke algebras and
quotients of them:
\DefLemma{P3.3b} Let $\cal A$ be a Hecke algebra as in \Ref{P3.3a}.
Let $\pi$ be an
element of the symmetric group $\Pi_n$ of degree~$n$, i.~e.~a permutation
of the set $I:=\{1, 2, \ldots, n\}$. Then~$\pi$ can
be written as the composition of transpositions $t_j$ (where $t_j$ interchanges
the elements $j$ and $j+1$ of $I$). The minimal number of such transpositions
is called {\sl length of $\pi$} and is denoted by $l(\pi)$.
Let $\pi=t_{k_1}\cdots t_{k_l}$ be a decomposition of $\pi$ into a
minimal number $l=l(\pi)$ of transpositions. Then $\sigma_{k_1}\cdots
\sigma_{k_l}$ does not depend on the actual choice of transpositions as
far as their number is minimal. Therefore
$\sigma_\pi:= \sigma_{k_1}\cdots \sigma_{k_l}$
is well-defined.
%(Remark: The elements $\sigma_\pi$ for $\pi\in \Pi_n$ form a linear
%basis of~$\cal A$.)
\Ende

\DefLemma{P3.3c} Let $\cal A$ be a Hecke algebra as in \Ref{P3.3a}.
Then define the element
$$S_n:= \sum_{\pi\in \Pi_n} q^{-2l(\pi)}\sigma_\pi\in {\cal A}.$$
This element satisfies the property $(\sigma_k-1)S_n = 0$ for $1\le k\le n-1$.
\Bew Let $k$ be an integer between~1 and $n-1$.
Let $\pi\in\Pi_n$ be a permutation such that $\pi(k)>\pi(k+1)$ and let
$\pi':=t_k\pi$. If $t_{k_1}\cdots t_{k_l}$
is a decomposition of~$\pi'$ into a minimal number of transpositions then
$t_kt_{k_1}\cdots t_{k_l}=t_k\pi'$ is a decomposition of~$\pi$ into a minimal
number of transpositions and $l(\pi)=l(\pi')+1$. Therefore for all $k$
{\def\Summe#1{\sum_{\textstyle{\pi\in\Pi_n \atop \!\!\pi(k)#1\pi(k+1)\!\!}}
q^{-2l(\pi)}\sigma_\pi}
$$\eqalignno{S_n&=\big(\Summe > +\Summe < \big) =\cr
&=\big(\sigma_kq^{-2} \Summe < +\Summe <\big) =
(q^{-2}\sigma_k+1)\Summe <}$$
and hence
$$(\sigma_k-1)S_n =
\underbrace{(\sigma_k-1)(1+q^{-2}\sigma_k)}_{=0 {\rm \ (
Hecke\ algebra)}}\Summe < =0.$$}
\Ende

\Rema{P3.3d} The Hecke algebra is a generalization of the symmetric
group, and for $q=1$ the Hecke algebra relations are just the relations
between the transpositions of the symmetric group. Let~$V$ be a
vector space. The symmetric group acts on~$V^{\tensor n}$ by
permutations of the tensor factors. The operator~$\sigma_k$ corresponding to
a transposition~$t_k$ has the eigenvalues~1 and -1. The intersection
of the kernels of all $\sigma_k-1$ or of the kernels of all $\sigma_k+1$
are called \Quote{totally symmetric vectors} or \Quote {totally
antisymmetric vectors}, respectively. When a Hecke algebra (or a quotient
of it) acts on~$V^{\tensor n}$,
then the eigenvalues are~1 and~$-q^2$ due to the second Hecke algebra
relation. The intersection of the kernels of all $\sigma_k-1$ or of
the kernels of all $\sigma_k+q^2$ is called the
space of \Quote{totally $q$-symmetric vectors} or \Quote{totally
$q$-antisymmetric vectors}. The element~$S_n$ is called
\Quote{symmetrization operator}, which is justified by \Ref{P3.3c},
which also explains
the factor $q^{-2l(\pi)}$ in the definition of~$S_n$.
\Ende

\parindent=0pt
Proof of \Ref{P3.2}: Let $K={}^2\C$ be the space of column vectors and $K'=\C^2$
the dual space of row vectors.
$w^0$ has always the matrix element~1 (because $\Delta(1)=1\tensor 1$).
Let $w$ be the fundamental \cor~$w^{1/2}$. Then
$w\tensor w\cong w^0\oplus w^1$, which is equivalent to the matrix
$$\pmatrix{1&0\cr 0&w'\cr},$$
where $w'\in M_3(\C)$.
This matrix has the column eigenvector $E=(1\ 0\ 0\ 0)^T$ and the
row eigenvector $E'= (1\ 0\ 0\ 0)$. Therefore there are the relations
$$E'(w\tensor w)=E' = w^0E', \ (w\tensor w)E=E = Ew^0.$$
Thus the vectors $E$ and $E'$, considered as $4\times 1$ oder $1\times 4$
matrices, intertwine $w\tensor w$ and $w^0$. Moreover
$$E'E\ne 0.  \Eqno 4$$
Now $(E'\tensor \Eins_2)(\Eins_2\tensor E)$ can be regarded as an
intertwiner of~$w$ with~$w$,
because $w\cong w\tensor w^0 \cong w^0\tensor w$ and
$$(E'\tensor \Eins_2)(\Eins_2\tensor E)(w\tensor w^0) =
(E'\tensor \Eins_2)(w\tensor w\tensor w)(\Eins_2\tensor E) =
(w^0\tensor w)(E'\tensor\Eins_2)(\Eins_2\tensor E).$$
Since $w$ is irreducible, by Schur's Lemma (\SRef{P2.14})
$(E'\tensor \Eins_2)(\Eins_2\tensor E)$ is a multiple of the identity,
say~$\lambda$ times the identity.
Using the coordinate representation of $E$ and $E'$ with respect to a
basis $\{e_i\tensor e_j\mid 1\le i,j\le 2\}$ of
$K\tensor K$ and the dual basis
$\{e_i'\tensor e_j'\mid 1\le i,j\le 2\}$ of $K'\tensor K'$,
$$E= \sum_{i,j} E_{ij}e_i\tensor e_j, \
E' = \sum_{i,j} E'_{ij}e_i'\tensor e_j',$$
this condition becomes
$$\sum_{k=1}^2 E'_{ik}E_{kj} = \lambda\delta_{ij}.$$
Therefore the matrices ~$\widehat E$ with entries~$E_{ij}$
and~$\widehat E'$ with entries $E'_{ij}$ satisfy
$$\widehat E'\widehat E = \lambda \Eins_2.$$
If $\lambda=0$ then $\widehat E$ must have rank~1, because if it
has rank 2, then $E'=0$ and if it has rank~0
then $E=0$ in contradiction to~(4). Hence $\widehat E_{ij}=x_iy_j$ for
some $x_1, x_2, y_1, y_2\in\C$
and~$E$ has the form
$E=x\tensor y$, where $x=x_1e_1+x_2e_2$, $y=y_1e_1+y_2e_2$. From $(w\tensor w)E=E$ it
follows that 
$$wx\tensor wy = x\tensor y, \ x\tensor wy=w\me x\tensor y.$$
Both sides are in~$K\tensor K\tensor A$. Applying $\phi\tensor\id_K\tensor \id_A$,
where $\phi$ is a
linear form on $K$ such that $\phi(x)=1$, we get
$$wy = y\tensor (\phi\tensor \id_A)(w\me x).$$
Therefore $\C y$ is an $w$-invariant subspace in contradiction to the
fact that~$w$ is irreducible.
Thus $\lambda\ne 0$, and by scaling of $E'$ which does not
change the relations, one gets $\widehat E'= \widehat E\me$.
The vector~$E$ in $K\tensor K$ can be written as sum of a symmetric
tensor~$\Esy$, i. e. an element of $K\tensor K$
which is invariant with respect
to the flip automorphism $\tau$ of $ K\tensor K$, mapping
$x\tensor y$ to $y\tensor x$, and an antisymmetric tensor~$\Easy$ satisfying
$\tau(\Easy)=-\Easy$, defined by $\Esy=\halb(E+\tau(E))$ and
$\Easy=\halb(E-\tau(E))$. Symmetric tensors $\sum_{i,j}a_{ij}e_i\tensor e_j$
in $K\tensor K$, where $a_{ij}=a_{ji}$ for all~$i$, $j$, can be identified with
quadratic forms~$Q$ on~$K'$, namely $Q(\sum_i v_ie_i') = \sum_{i,j} a_{ij}v_iv_j$.
In particular there are  bases such that~$\Esy$ has one of the following
presentations:
\item (a) $\Esy= e_1\tensor e_2+e_2\tensor e_1$ if $Q$ has rank 2,
\item (b) $\Esy = e_1\tensor e_1$ if $Q$ has rank 1,
\item (c) $\Esy = 0$ if $Q$ has rank 0.

With respect to any basis $\{\tilde e_1, \tilde e_2\}$ of $K$ an antisymmetric
tensor~$\Easy$ is a scalar multiple of $\tilde e_1\tensor \tilde e_2-\tilde e_2
\tensor \tilde e_1$.
Therefore~$E$ has one of the following presentations:
\item (a) $E= (1+c)e_1\tensor e_2+(1-c)e_2\tensor e_1$ with $c\in \C$.
Since $\widehat E$ has rank~2, both coefficients must not vanish. Therefore~$E$
is a scalar
multiple of $e_1\tensor e_2-qe_2\tensor e_1$, where $q={c-1\over c+1}$
and $q\in \C\setminus\{0,1\}$.
\item (b) $E= e_1\tensor e_1+c(e_1\tensor e_2-e_2\tensor e_1)$, where
$c\in\C\setminus\{0\}$, because $\widehat E$ has rank~2. Therefore~$E$ is
a scalar multiple of $e_1\tensor e_2-e_2\tensor e_1+t e_1\tensor e_1$,
where $t={1\over c}$. According to \Ref{P3.3} this is equivalent to
the vector for $SL_{t=1}(2)$. In this case let $q:=1$.
\item (c) $E= c(e_1\tensor e_2-e_2\tensor e_1)$, where $c\in \C\setminus\{0\}$.
This is the case $q=1$ which is not included in~(a).

Now let the associative, unitary algebra $\Ho$ be generated by
the elements $\alpha$, $\beta$, $\gamma$, $\delta$ subject
to the relations $(v\tensor v)E=E$, $E'(v\tensor v)=E'$ for
$v=({\alpha\atop\gamma}\ {\beta\atop\delta})$.
There are uniquely determined comultiplication, counit and antipode such that
this algebra becomes a Hopf algebra and~$v$ is a \cor
(see \Ref{P4.2} and \Ref{P4.3} below).

Since the relations between the generators of $\Ho$ are satisfied
in $\cal H$, there is a Hopf algebra map
$$\psi\colon {\Ho}\to {\cal H}, v_{ij}\mapsto w_{ij}.$$
We shall study the \cor theory of~$\Ho$.

Consider the $4\times 4$ matrix $\sigma:=\Eins_4 +qE\cdot E'$ (where
$E$ and~$E'$ are again $4\times 1$ and $1\times 4$ matrices, respectively). Then
$\sigma$ is an element of the vector space
$\Mor(v\tensor v, v\tensor v)$. It satisfies the
relations
$$(\sigma-\Eins_4)(\sigma+q^2\Eins_4) = 0, \
(\sigma\tensor\Eins_2)(\Eins_2\tensor \sigma)(\sigma\tensor \Eins_2)=
(\Eins_2\tensor \sigma)(\sigma\tensor \Eins_2)(\Eins_2\tensor \sigma).$$
Fix an integer~$n\ge 2$ and define for integers~$k$ satisfying $0<k<n$:
$$\sigma_k = \underbrace{\Eins_2\tensor\cdots\tensor\Eins_2}_{k-1}
\tensor \sigma \tensor \underbrace{\Eins_2\tensor\cdots\tensor\Eins_2}_{n-k-1}.$$
These are operators on the $n$-fold tensor product $K^{\tensor n}$
and intertwine~$v^{\tensor n}$ with $v^{\tensor n}$.
They satisfy the {\sl Hecke algebra relations} (cf. \Ref{P3.3a}).

Now define the operators $\sigma_\pi$ as in \Ref{P3.3b} and
the {\sl symmetrization operator} as in \Ref{P3.3c}:
$$S_n:= \sum_{\pi\in \Pi_n} q^{-2l(\pi)}\sigma_\pi.$$
Due to \Ref{P3.3c} it takes values in
$$K^{n/2}:= \{x\in K^{\tensor n}\mid
\forall k\colon \sigma_k(x) = x\} = \bigcap_{k}\Ke(\sigma_k-\Eins).$$

The dimension of the space $K^{n/2}$ is $n+1$.
(Proof: analyze relations on coordinates of elements of $K^{n/2}$ or
see \Wii{}).
The space $K^{n/2}$ is $v^{\tensor n}$-invariant as intersection of the
kernels of the intertwiners $\sigma_k-\Eins$ by \Ref{P2.14n} and \Ref{P2.14o}
and a right comodule. Let~$v^{n/2}$ denote the corresponding sub\cor of~$v^{\tensor n}$ as in \Ref{P2.14m}.
Then $v^{n/2}$ is a \cor of dimension $n+1$.
By definition, $v=v^{1/2}$, $v^0$ is the trivial one-dimensional
\cor. At this moment we assume that~$q$ is not a non-real root of unity.
For all $s\in \halb \N_0$, one has
\item (A) $v^k$ is irreducible for all $k\le s+\halb$,
\item (B) $v^k\tensor v\cong v^{k+\halb}\oplus
v^{k-\halb}$, where by definition $v^{-1/2}:=0$, for all $k\le s$.

These statements will be proved by induction:
The case $s=0$ follows from the result on monomials below.

Suppose the statements are true for $s$ replaced by $s-\halb$.

We want to decompose $v^s\tensor v$ and consider the map
$$\phi\colon K^{s-\halb}\to K^s\tensor K,
x\mapsto (S_{2s}\tensor \Eins_2)(x\tensor E)$$
(antisymmetrization-symmetrization procedure).
Note that $E\in K\tensor K$. The map is well defined due to the property
of the symmetrization operator and intertwines $v^{\tensor(2s-1)}\cong
v^{\tensor (2s-1)}\tensor v^0$ with $v^{\tensor 2s}\tensor v$.
By inspection
$$(\sigma_{2s}-\Eins)\phi\underbrace{(e_1\tensor
\cdots\tensor e_1)}_{2s-1 {\rm \ factors}}\ne 0$$
if $q$ is not a non-real root of unity. Therefore $\phi(e_1\tensor\cdots \tensor e_1)\notin
K^{s+\halb}$ and $\Ke(\phi)\ne K^{s-\halb}$.
By induction hypothesis, there is no proper non-trivial $v^{2s-1}$-invariant
subspace of $K^{s-\halb}$ and
the kernel of $\phi$ is invariant by \Ref{P2.14o}. Consequently~$\phi$ is
injective and $\Im(\phi)$ corresponds to~$v^{s-\halb}$.
Moreover
$$\Im(\phi)\cap K^{s+\halb}\neq \Im(\phi).$$
Since $v^{s-\halb}$ is irreducible, there is no proper non-trivial $v^{s-\halb}$-invariant
subspace of $\Im(\phi)$ and
$\Im(\phi)\cap K^{s+\halb}=\{0\}$. Thus
$$K^{s-\halb}\oplus K^{s+\halb} \cong \Im(\phi)\oplus K^{s+\halb}
\subset K^s\tensor K.$$
Equality follows by dimension arguments ($\dim K^t = 2t+1$) and yields~(B).
By definition of the tensor product of representations, the monomials
in $\alpha$, $\beta$, $\gamma$, $\delta$ of
degree smaller or equal to $2s+1$ are linear combinations of the matrix
elements of $v^{\tensor(2s+1)}$. Using result (B) yields that they are
linear combinations of matrix elements of $v^0, v^\halb, \ldots, v^{s+\halb}$.
The space of these monomials has dimension $\sum_{k=1}^{2s+2} k^2$,
as in the classical case for~${\rm Pol}(SL(2, \C))$.
This has been shown in \Wii{}, \WZ{}.
Since $v^t$ has $(2t+1)^2$ matrix elements for all $t\in \halb \N_0$,
the space spanned by the matrix elements of
$v^0, v^\halb, \ldots, v^{s+1/2}$ has this dimension if and only if all matrix
elements are linearly independent. Hence (A) follows. Now it is
easy to prove that $\Ho$ is a quantum $SL(2)$-group.

The matrix $\psi(v^s)$ is a \cor of $\cal H$ because~$\psi$
respects $\Delta$ and $\eps$.

Then $\psi(v^s)\cong w^s$ for
all $s\in \halb\N_0$.
Proof by induction:
The assertion is trivial for $s=0$ and $s=\halb$.
Suppose the statement is true for all non-negative half integers smaller
than~$s$.
Then by part (B), $v^{s-\halb}\tensor v\cong v^{s-1}\oplus v^{s}$.
Since $\psi$ is an algebra homomorphism, by the definitions of direct
sum and tensor product of \cors the following holds
$$\eqalign{w^{s-1}\oplus \psi(v^s)&\cong\psi(v^{s-1})\oplus\psi(v^s)\cong \psi(v^{s-1}\oplus v^s)\cong
\psi(v^{s-\halb}\tensor v) =\cr 
&= \psi(v^{s-\halb})\tensor \psi(v) \cong
w^{s-\halb}\tensor w \cong w^{s-1}\oplus w^s
.}$$
Due to condition (4) for the quantum $SL(2)$-group, the \cor
$\psi(v^s)$ is completely reducible, whence by \Ref{P2.16} $\psi(v^s)\cong w^s$.
Thus~$\psi$ is an isomorphism and $\cal H$ can be identified with~$\Ho$.

Now we consider the case when~$q$
is a non-real root of unity (see e.~g.~\KP{}).
Let $q$ be a non-real root of unity of order~$N$. Define
$$N_0:=\cases{N& if $N$ is odd, \cr N/2& if $N$ is even.}$$
Then $\Ho$ has a \cor\empty\foot{cf. \Tak{part 5.2}}
$$z=\pmatrix{\alpha^{N_0}&\beta^{N_0}\cr \gamma^{N_0}&
\delta^{N_0}}.$$
Then $v^k\tensor v\cong v^{k+\halb}\oplus v^{k-\halb}$ for
$k<(N_0-1)/2$, $v^k$ is irreducible for $k\le \halb(N_0-1)$, and
$$v^{(N_0-1)/2}\tensor v \cong
\pmatrix{ v^{\halb N_0-1}&*&*\cr 0&z&*\cr 0&0&v^{(N_0-1)/2}}.$$
It is possible to show $\psi(v^k)\cong w^k$ for $k\le \halb(N_0-1)$ as before, but
on the other hand
$$\eqalign{w^{\halb N_0-1}\oplus w^{\halb N_0}&\cong w^{\halb(N_0-1)}\tensor w \cong
\psi(v^{\halb(N_0-1)}\tensor v) \cong\cr
&\cong
\pmatrix{\!w^{\halb N_0-1}\!&*&*\cr 0&\!\psi(z)\!&*\cr 0&0&w^{\halb N_0-1}}
\cong w^{\halb N_0-1}\oplus \psi(z)\oplus
\!w^{\halb N_0-1}\!,}$$
because \cors in $\cal H$ are completely reducible. But this is
a contradiction to \Ref{P2.16}.

Let $q_1$ and $q_2$ be two values such that $SL_{q_1}(2)\cong SL_{q_2}(2)$
($q_1, q_2\in\C\setminus Y \cup \{t=1\}$, where the subset~$Y$ contains~0 and
all non-real roots of unity),
i.~e.~that the Hopf algebras are isomorphic. Then the fundamental
representation $w_1$ is mapped to $w_2$, i.~e.~they are equivalent:
$w_1 = Qw_2Q\me$. Let $E_1$, $E_2$ be the corresponding eigenvectors.
Then
$$\displaylines{(w_1\tensor w_1)E_1= E_1, (Qw_2Q\me\tensor Qw_2Q\me)E_1=E_1 \Folgt\cr
(w_2\tensor w_2)((Q\me\tensor Q\me)E_1)=(Q\me\tensor Q\me)E_1  = \lambda E_2}$$
with $\lambda\in\C\setminus\{0\}$,
because the space of eigenvectors of $w_2\tensor w_2$ for the eigenvalue~1
is one-dimensional.
There are symmetric tensors $\Esy_1$, $\Esy_2$ and antisymmetric
tensors $\Easy_1$, $\Easy_2$ such that $E_1= \Esy_1+\Easy_1$
and $E_2=\Esy_2+\Easy_2$. Therefore
$$E_1=\lambda(Q\tensor Q)E_2 \Folgt \Esy_1=\lambda(Q\tensor Q)\Esy_2,
\Easy_1=\lambda (Q\tensor Q)\Easy_2\Eqno 5$$
and $\Esy_1$ and $\Esy_2$ have the same rank. If the rank is~0 or~1,
it is the same deformation, and if the rank is~2, one can use (\SRefg{5}) and
the fact that the rank of $(Qe_1\tensor Qe_2)$ is one, to
get $q_1=q_2$ or $q_1q_2=1$ (in the last case the isomorphism is
given by $e_1\leftrightarrow e_2$).

\Abschnitt{Quantum $SL(N)$-groups}
Let $N$ be a positive integer greater than 1.
The Hopf algebra $\cal H$ of the group $SL(N,\C)$ corresponds to
the commutative unital
algebra generated by the matrix elements~$w_{ij}$ for $1\le i,j\le N$
of a fundamental \cor~$w$ subject to the relations
$$w^{\tensor N}E=E, \ E'w^{\tensor N}=E' \Eqno {G5s}$$
where $E$ and $E'$ are classical completely antisymmetric elements
of $(^N\C)^{\tensor N}$ and $(\C^N)^{\tensor N}$ respectively,
i.~e.~with respect to a basis $\{e_1, \ldots, e_N\}$ of~$^N\C$ and
a dual basis $\{e_1', \ldots, e_N'\}$ of~$\C^N$, they can be presented
as
$$E=\sum_{\pi\in \Pi_n} (-1)^{l(\pi)} e_{\pi(1)}\tensor
\cdots \tensor e_{\pi(N)}, \
E'=\sum_{\pi\in \Pi_n} (-1)^{l(\pi)} e'_{\pi(1)}\tensor
\cdots \tensor e_{\pi(N)}'. \Eqno{G5ss}$$
Then the relations just mean 
(assuming commutativity) that the determinant of the matrix~$w$ is one.
For $SL(2)$ this is just $e_1\tensor e_2-e_2\tensor e_1$ and
$e_1'\tensor e_2'-e_2'\tensor e_1'$, which is changed to
$e_1\tensor e_2-qe_2\tensor e_1$ and up to a
non-zero factor to $e_1'\tensor e_2'-q e_2'\tensor e_1'$
in the standard deformation $SL_q(2)$. Therefore it is natural to define
$$E_q=\sum_{\pi\in \Pi_n} (-q)^{l(\pi)} e_{\pi(1)}\tensor
\cdots \tensor e_{\pi(N)}, \
E'_q=\sum_{\pi\in \Pi_n} (-q)^{l(\pi)} e'_{\pi(1)}\tensor
\cdots \tensor e_{\pi(N)}' \Eqno{G5ta}$$
and to consider the relations
$$w^{\tensor N}E_q=E_q, \ E'_qw^{\tensor N}=E'_q. \Eqno {G5tb}$$
For $q$ not being a non-real root of unity they imply (cf. \TK{})
$$ w^{\tensor 2}\sigma=\sigma w^{\tensor 2} \Eqno{G5tc}$$
where
$$\sigma(e_i\tensor e_j):=\cases{qe_j\tensor e_i &if  $ i<j$,\cr
qe_j\tensor e_i+(1-q^2)e_i\tensor e_j&if $i>j$,\cr
e_i\tensor e_i&if $i=j$,}$$ 
for $i,j=1,\dots,N$. Now $SL_q(N)$ is introduced as the unital
algebra generated by $w_{ij}$ for $1\leq i,j\leq N$ subject to
the relations $(\SRefg{G5tb})$, $(\SRefg{G5tc})$ (cf. \Pdiii{}).
One can check that this definition coincides with the standard
one (cf. \Drin{}, \Ross{}).

The following proposition shows that all unital algebras with
relations defined by intertwiners are bialgebras.
If the intertwiners are chosen badly, the bialgebras can be
small and uninteresting. For each matrix~$w$ and each~$n\in \N$
define the matrix~$w^{\tensor n}$ as for \cors in \Ref{P2.13}
and let $w^{\tensor 0}:=\Eins_1$.
\Prop{P4.2}
Let $\cal H$ be the universal unital algebra generated by elements~$w_{ij}$
for $1\le i,j\le N$, which are the entries of a matrix~$w$ subject
to relations
$$%E_m(\Eins_N^{\tensor t_m}\tensor w^{\tensor s_m})=(w^{\tensor t_m}
%\tensor \Eins_N^{\tensor s_m})E_m \hbox{ or shortly }
E_m w^{\tensor s_m} = w^{\tensor t_m} E_m \Eqno 6$$
for $m$ in an index set~$I$, $s_m,t_m\in\N_0$ and
$E_m\in %({}^N\C)^{\tensor t_m}\tensor (\C^N)^{\tensor s_m}
M_{N^{t_m}\times N^{s_m}}(\C)$.
Then there exist a unique comultiplication
and counit such that~$\cal H$ is a bialgebra and~$w$ is a \cor
of~$\cal H$.
\Bew \voritem (a) Uniqueness: We must have $\Delta w_{ij}=\sum_{k=1}^N w_{ik}\tensor w_{kj}$
and $\eps(w_{ij})=\delta_{ij}$ for all $i$ and~$j$. Since~$\Delta$
and~$\eps$ are unital algebra homomorphisms, they are uniquely determined if they exist.
\item (b) Existence: Define $\widehat w_{ij}:= \sum_{k=1}^N w_{ik}\tensor_{\Bbb C}
w_{kj}\in {\cal H}\tensor {\cal H}$ for all~$i$ and~$j$. The matrix $\widehat w$ with entries
$\widehat w_{ij}$ also satisfies the relations (\SRefg{6}), because
$\widehat w^{\tensor n} = w^{\tensor n}\tensor_{\Bbb C} w^{\tensor n}$
follows from the rule $(a\tensor b)(c\tensor d)=(ac\tensor bd)$ and
$$\eqalign{E_m\widehat w^{\tensor s_m} &= E_m(w^{\tensor s_m} \tensor_{\Bbb C} w^{\tensor s_m})=
w^{\tensor t_m} E_m \tensor_{\Bbb C} w^{\tensor s_m} =\cr
&=w^{\tensor t_m} \tensor_{\Bbb C} E_m w^{\tensor s_m} =
w^{\tensor t_m}\tensor_{\Bbb C} w^{\tensor t_m} E_m =
\widehat w^{\tensor t_m} E_m,}$$
because the entries of $E_m$ are just complex numbers.
Define $\widetilde w_{ij}:= \delta_{ij}$ for all $i$, $j$. Then the
matrix $\widetilde w$ with entries $\widetilde w_{ij}$ satisfies the
properties
$$E_m \widetilde w^{\tensor s_m} = E_m, \
  \widetilde w^{\tensor t_m}E_m = E_m,$$
whence it satisfies relations (\SRefg{6}). Now the universality of~$\cal H$
gives the existence of unital homomorphisms~$\Delta$, $\eps$ such
that $\Delta(w_{ij})=\widehat w_{ij}$ and $\eps(w_{ij})=\widetilde w_{ij}$.
It is enough to check Conditions~(1) and~(2) for bialgebras (cf. \Ref{P2.6})
for elements $f=w_{ij}$ when they are obvious.
\Ende

\Prop{P4.3} Let the conditions of \Ref{P4.2} be satisfied. Let
$\{e_1, \ldots, e_N\}$ be a basis of~$^N\C$ and $\{e_1', \ldots,
e_N'\}$ be a dual basis of~$\C^N$. Moreover
assume that there exist positive integers~$s$ and~$t$ and
elements $E\in \Mor(\Eins_1 , w^{\tensor t})$
and $E'\in\Mor( w^{\tensor s}, \Eins_1)$ such that
$$E= \sum_{k=1}^N e_k\tensor f_k, \
E'= \sum_{k=1}^N f_k'\tensor e_k'$$
such that the elements $f_k\in (^N\C)^{\tensor t-1}$ and
$f_k'\in (\C^N)^{\tensor s-1}$ are linearly independent.
Then the matrix~$w\me$ exists and there is a uniquely determined antipode~$S$
such that the bialgebra $\cal H$ is a Hopf algebra.
\Bew From the relation $w^{\tensor t}E=E$ it follows that
$$(w\tensor w^{\tensor (t-1)})E=E \Folgt
\sum_{k=1}^N we_k\tensor w^{\tensor (t-1)}f_k = \sum_{k=1}^N e_k\tensor f_k.\Eqno{G4.3}$$
Since the elements~$f_k$ of $(^N\C)^{\tensor (t-1)}$
are linearly independent, there are elements~$g_k'$ of the dual space
$(\C^N)^{\tensor (t-1)}$ such that $g_i'f_j=\delta_{ij}$.
Apply $e_i'\tensor g_j'$ to \Refg{G4.3}:
$$\sum_{k=1}^N e_i'we_k\tensor g_j'{w^{\tensor (t-1)}{f_k}} =
\sum_{k=1}^N e_i'e_k\tensor g_j'f_k= 1\tensor g_j'f_i = \delta_{ij}.$$
Therefore the matrix~$G$ with entries $G_{kj}:=g_j'w^{\tensor (t-1)}f_k$ is a
right inverse to~$w$.
{}From the second condition it follows in a similar way that there is a left
inverse of~$w$. Thus $w\me$ exists.
Finally, when to the relation
$$E_m w^{\tensor s_m} = w^{\tensor t_m} E_m ,$$
$(w^{\tensor s_m})\me = (w\me)^{\tensorop s_m}$ is applied to the right and
$(w^{\tensor t_m})\me = (w\me)^{\tensorop t_m}$ to the  left (the
tensor product \Quote{$\tensorop$} is $\tensor$ with respect to
the algebra ${\cal H}^{\rm op}$ with opposite multiplication), then
$$(w\me)^{\tensorop t_m}E_m=E_m (w\me)^{\tensorop s_m}.$$
Therefore
there is a unital algebra homomorphism~$S\colon {\cal H}\to {\cal H}^{\rm op}$
such that $S(w)=w\me$. Equivalently, $S\colon {\cal H}\to {\cal H}$
is a unital antihomomorphism. It is enough to check Condition~(3) for
the antipode (cf. \Ref{P2.6}) for $f=w_{ij}$ when it is obvious.
Uniqueness of~$S$ follows from \Ref{P2.7}.
\Ende

\Rema{P4.4} For the quantum $SL(N)$ group take $I=\{1, 2, 3\}$,
$E_1=E_q$, $t_1=N$, $s_1=0$, $E_2=E'_q$, $t_2=0$, $s_2=N$,
$E_3=\sigma$, $t_3=s_3=2$. Then the algebras
$SL_q(N)$ are Hopf algebras.
\Ende

\Rema{P4.1} For $0<q\le 1$ the corepresentation theory of $SL_q(N)$
is the same as for the classical $SL(N)$ (cf.
\TK{}, \Pdiii{}).
If~$q$ is transcendental, see \Ross{}, \Hay{}.
If $q\in \C\setminus\{0\}$ is not a non-real root of unity, see
\PW{}.
There are deformations of the
orthogonal and symplectic groups \FRT{}, \Ta{} (cf. \Pdiii{}).
\Ende

\subsection{P5}{$*$-Structures}
In the classical theory there exist $*$-structures on $\Pol(SL(2))$ which
give  the Hopf \hbox{$*$-algebras} $\Poly(SU(2))$, $\Poly(SU(1,1))$ and $\Poly(SL(2, \R))$.
We will classify the Hopf \hbox{$*$-algebra} structures on the quantum $SL(2)$-groups~$\cal H$
described in \Ref{P3.2}. Firstly recall that~$\cal H$ is generated as an algebra
by the matrix elements of a $2\times 2$ matrix~$w$ subject to the
relations
$$(w\tensor w)E=E, \ E'(w\tensor w)=E'$$
or equivalently
$$\sum_{j,l} w_{ij}w_{kl}E_{jl} = E_{ik}, \
\sum_{i,k} E_{ik}'w_{ij}w_{kl} = E_{jl}'. \Eqno 7$$

\Lemma{P5.1o}
Let $\psi$ be an (anti-)linear comultiplicative algebra
(anti-)automorphism of a quantum $SL(2)$-group ${\cal H}$. Then
\item (a) there exists a matrix $Q\in GL(2,\C)$ such that $\psi(w)=Q wQ\me$.
\item (b) If and only if the matrix $Q\in GL(2,\C)$ satisfies the conditions
$$(Q\me\tensor Q\me)E=cE, \ E'(Q\tensor Q)=c'E' \Eqno{G7.2}$$
for some numbers $c,c'\in \C\setminus\{0\}$,
there is a Hopf algebra automorphism~$\psi$ of~$\cal H$ such that $\psi(w)=Q wQ\me$.
Moreover, all Hopf algebra automorphisms of~$\cal H$ can be described in this way.
\item (c) Let $\tau$ denote the linear twist
(interchanging factors) and let~$\bar E$ and~$\bar E'$ denote the elements
of $^2\C\tensor {}^2\C$ and $\C^2\tensor \C^2$ with conjugate complex
coefficients with respect to the bases $e_i\tensor e_j$, $e_i'\tensor e_j'$. Then if and only if the matrix $Q\in GL(2,\C)$ satisfies the conditions
$$(Q\me\tensor Q\me)\tau \bar E=c E, \ \bar E'\tau
(Q\tensor Q)=c' E' \Eqno{G7.2a}$$
for some $c,c'\in \C\setminus\{0\}$,
there is an antilinear, comultiplicative, algebra antiautomorphism~$\psi$ of~$\cal H$
such that $\psi(w)=Q wQ\me$.
\item (d) Let the antilinear involutive
comultiplicative algebra antiautomorphisms~$\psi$, $\hat\psi$ and the
corresponding matrices~$Q$ and~$\hat Q$ be defined as in (a).
Then the Hopf algebra $\cal H$ equipped with $*$-structures $\psi$ and $\hat\psi$ 
gives isomorphic Hopf $*$-algebras if and
only if $\hat\psi$ is
equivalent to $\psi$ up to
a Hopf algebra automorphism~$\phi$ (i.~e.~$\hat\psi=\phi\psi\phi\me$)
if and only if
$\hat Q=c\bar A\me QA$ where $c\in \C\setminus\{0\}$ and
$A\in GL(2, \C)$ corresponds to~$\phi$ via (b).
\Bew \voritem (a) Since $\psi$ is comultiplicative,
the matrix $\psi( w)$ is a \cor. The following conclusions
follow from the fact that $\psi$ is bijective:
$w$ is irreducible if and only if the matrix elements~$w_{ij}$
are linearly independent if and only if the matrix elements~$\psi(w_{ij})$
are linearly independent if and only if~$\psi( w)$ is irreducible.
But there is only one irreducible \cor of dimension~2 up
to isomorphism, therefore there is a matrix~$Q\in GL(2,\C)$ such that
$$\psi( w) = QwQ\me. \Eqno{7.1}$$
\item (b) Since the trivial \cor appears in the direct sum
decomposition $w\tensor w\cong w^1\oplus w^0$ only once, by \Ref{P2.14}
the space of intertwiners in $\Mor(w\tensor w,w^0)$ is one-dimensional.
Thus Condition (\SRefg{G7.2}) is equivalent to the condition
that ${(Q\me\tensor Q\me) E}$ intertwines $w\tensor w$ with~$w^0$ and
$E'(Q\tensor Q)$ intertwines $w^0$ with $w\tensor w$:
$$\left.\eqalign{(w\tensor w)(Q\me\tensor Q\me) E &=(Q\me\tensor Q\me) E \iff
(QwQ\me\tensor QwQ\me)E = E \cr
\hbox{ and }E'(Q\tensor Q)(w\tensor w)&=E'(Q\tensor Q) \iff
E'(QwQ\me\tensor QwQ\me)=E'.}\right\} \Eqno{G7.3}$$
Let~$\psi$ be a Hopf algebra automorphism of~$\cal H$. Then by part~(a),
there is a matrix $Q\in GL(2,\C)$ such that $\psi(w)= Qw Q\me$.
The automorphism~$\psi$ must map the relations between the generators
of~$\cal H$ to relations in~$\cal H$, therefore \Refg{G7.3} holds.
Conversely, let \Refg{G7.3} be satisfied.
Let~$\cal F$ be the
free associative unital algebra generated by the matrix elements of~$w$
and let~$\cal I$ be the two-sided ideal generated by the relations
(\SRefg{7}). Then the map $\psi$ can be defined as unital algebra
homomorphism on~$\cal F$ such that $\psi(w)=QwQ\me$.
Equation (\SRefg{G7.3}) shows
that~$\psi$ maps $\cal I$ to $\cal I$, therefore it induces a unital
algebra homomorphism on $\cal H= F/I$. Such~$\psi$ preserves the
Hopf algebra structure of~$\cal H$. Moreover, replacing~$Q$ by~$Q\me$
(Equation (\SRefg{G7.2}) still holds for $c\me$ and $(c')\me$)  we get $\psi\me$.
\item (c) The proof is similar as in part (b). The only changes arise from
the fact that~$\psi$ should be an antilinear algebra antiautomorphism instead
of a linear algebra automorphism. Therefore, $\psi$ applied to relations
(\SRefg{7}) yields
$$\sum_{j,l} \psi(w_{kl})\psi(w_{ij})\bar E_{jl} = \bar E_{ik}, \
\sum_{i,k} \bar E_{ik}'\psi(w_{kl})\psi(w_{ij}) = \bar E_{jl}'  $$
or shortly
$$\tau(\psi(w)\tensor \psi(w))\tau\bar E = \bar E, \
\bar E'\tau(\psi(w)\tensor \psi(w))\tau = \bar E'.$$
Using $\tau^2=\id_{\cal H}$, we get the desired results.
\item (d) $\phi\psi\phi\me(w)=\phi\psi(A\me w A)=\phi(\bar A\me QwQ\me \bar A)=
\bar A\me QAwA\me Q\me\bar A$, while
$\hat\psi(w)=\hat Q w\hat Q\me$. The left hand sides are equal if and
only if $\hat Q\me\bar A\me QA\in \Mor(w,w)=\C\Eins_2$.
\Ende
Remark. It is easy to check that the second condition in (\SRefg{G7.2})
(and also the second condition in (\SRefg{G7.2a})) is redundant.

\Thm{P5.1} All non-equivalent Hopf $*$-algebra
structures on
the quantum $SL_q(2)$-groups ${\cal H}$  are defined by
$\bar w = QwQ\me$, where
\item (a) $Q=\Smatrix{1\atop 0\cr 0\atop 1}$, $\vert q\vert=1$. Then $\bar w=w$.
This algebra is called
$\Poly(SL_q(2, \R))$.
\item (b) $Q= \Smatrix{0\atop 1\cr q\atop 0}$, $q\in \R\setminus\{0\}$.
Then
$w^* B w = wBw^*=B$,
for $B:=\Smatrix{1\atop 0\cr 0\atop -1}$.
This algebra is called
$\Poly(SU_q(1,1))$.
\item (c) $Q =\Smatrix{0\atop 1\cr -q\atop 0}$, $q\in\R\setminus\{0\}$.
Then
$w$ is unitary.
This algebra is called $\Poly(SU_q(2))$.

The only equivalence among them is $\Poly(SL_1(2, \R))\cong\Poly(SU_1(1,1))$.

For the non-standard deformation $SL_{t=1}(2)$ 
there is only one Hopf $*$-algebra structure (up to equivalence),
namely for $Q=\Smatrix{1\atop0\cr0\atop -1}$.
 
Except for (c), the above \cors $w$ are not equivalent to  unitary ones (The above examples
were given in \Wiii{}, \FRT{}, \Wii{}).
\smallskip
\rm Ideas of the proof: Since the map \Quote{$*$} is an antilinear comultiplicative
algebra antiautomorphism,
by \Ref{P5.1o}, part~(a) there is a matrix $Q\in GL(2,\C)$ such that
$\bar w = QwQ\me$.
By part~(c) of \Ref{P5.1o}, the map \Quote{$*$} can be an algebra
antiautomorphism if and only if~$Q$ satisfies the condition
$$(Q\me\tensor Q\me)
\tau \bar E =cE $$
for some $c\in\C\setminus\{0\}$.
The equation $*^2=\id_{\cal H}$ is equivalent to
$$\bar Q Q=d\Eins_2$$
with $d\in \C\setminus\{0\}$.
$Q$ is determined up to the equivalence relation as in \Ref{P5.1o}, part~(d).

Consider the standard quantum deformations $SL_q(2)$, $q\ne 1$, first.
{}From the relations (\SRefg{7}) it follows that
there are only the following
characters (algebra homomorphisms)~$\chi\colon {\cal H}\to \C$:
$$\chi_a(w)=\pmatrix{a&0\cr 0&a\me} \hbox{ and in addition to that for $q=-1$: }
 \chi_a'(w) = \pmatrix{0&a\cr a\me&0},$$
 where $a\in \C\setminus\{0\}$
(Relations (\SRefg{7}) are equivalent to
$$\displaylines{w_{11}w_{12}=qw_{12}w_{11}, \ w_{11}w_{21}=qw_{21}w_{11}, \
w_{12}w_{22}=qw_{22}w_{12}, \cr
w_{21}w_{22}=qw_{22}w_{21}, \
w_{12}w_{21}=w_{21}w_{12}, \cr
w_{11}w_{22}-qw_{12}w_{21}=
w_{22}w_{11}-q\me w_{12}w_{21}=1, }$$
and the numbers $\chi(w_{ij})$ should satisfy the same relations).

Now the following trick can be used in order to compute all possible
$*$-structures: If $\chi$ is a character, then also the map
$\chi^\#\colon x\mapsto \overline{\chi(x^* )}$ is a character,
because \C is commutative.

Then for any $a\in\C\setminus\{0\}$
there exists $b\in \C\setminus\{0\}$ such that $\chi_a^\#=\chi_b$ or
(for $q=-1$)
$\chi_a^\#=\chi_b'$. Applying both sides to~$w$, we get that~$Q$ is a diagonal
or antidiagonal matrix. Similarly (use $\chi\mapsto \chi\circ\phi$),
isomophisms~$\phi$ of Hopf algebras are given by diagonal or ($q=-1$)
antidiagonal matrices. Then we use the other conditions for~$Q$ and part~(d)
of \Ref{P5.1o}.

For the non-standard deformation $SL_{t=1}(2)$ split~$E$ into
$\Esy$ and $\Easy$ as in the proof of \Ref{P3.2}. Then consider~$Q$
with respect to both. For~$q=1$, equivalent~$Q$'s can be regarded as
matrices of the same antilinear mapping~$j$ such that~$j^2=d\cdot\id$
($j$ is equivalent to $kj$ for some $k\in \C\setminus\{0\}$).
Then~$d=1$ corresponds to~(a), (b) while $d=-1$ to~(c).
\Ende

\doppelpunktfalse \Rema{P5.4} \FRT{}, \Pdiii{}.
\item (a) There exist the following $*$-structures on $SL_q(N)$:
\itemitem (i) For $\vert q\vert =1$ you can choose $\bar w=w$.
The corresponding quantum group
is called $SL_q(N, \R)$.
\itemitem (ii) If~$q$ is real then
for $\eps_1, \ldots, \eps_N\in \{\pm 1\}$ there are $*$-structures
such that $w^*Bw=wBw^*=B$, where $B$ is a diagonal matrix with diagonal
elements $\eps_1, \ldots, \eps_N$. The corresponding quantum group
is called $SU_q(N; \eps_1, \ldots,\eps_N)$. For $\eps_1=\cdots=\eps_N=1$
we get the quantum group~$SU_q(N)$, in which~$w$ is a unitary
\cor.
\item (b) There are also $*$-structures on the orthogonal and symplectic
quantum groups. 
\Ende

\subsection{P6}{Compact Hopf $*$-algebras }

In this chapter we follow \Wi{}, \TK{}, \Koo{}.
Let $\cal A$ be a Hopf $*$-algebra.
\Def{P6.1} $\cal A$ is called {\sl compact} if there are unitary
\cors such that their matrix elements
generate $\cal A$ as algebra.
\Ende
Example: The fundamental \cor of $\Poly(SU_q(N))$ is unitary
and generates $\Pol(SL_q(N))$ as algebra.
\Lemma{P6.2} Let $\cal A$ be a compact Hopf $*$-algebra.
\item(a) The matrix elements of unitary \cors span $\cal A$.
\item (b) Let $v$ be a unitary \cor. Then $v$ is equivalent to
a direct
sum of irreducible unitary \cors.
\item (c) The matrix elements of non-equivalent irreducible unitary
\cors form a linear basis of $\cal A$.
\item (d) Each irreducible \cor is
equivalent to a unitary one.
\item (e) Each \cor is completely reducible (into irreducible
ones). Since the irreducible \cors are equivalent to unitary
\cors, each \cor is equivalent to a unitary
\cor.
\Bew \voritem (a) By definition, $\cal A$ is spanned by matrix elements
of tensor products of unitary \cors, but tensor products
of unitary \cors are unitary.
\item (b) Proof by induction with respect to the dimension~$d$ of \cors.
If $d=1$ or the \cor is irreducible, then there is nothing to
do.
Now assume that the \cor~$v$ is not irreducible. Then choose an
orthonormal basis of an invariant proper subspace~$L$ and add some more
orthonormal elements in order to get an orthonormal
basis~$B$ of $\C^{\dim v}$. The transition from the standard basis
to~$B$ is unitary and intertwines~$v$
with a unitary corepresentation
$$\pmatrix{A&B\cr 0&C\cr}=: w,$$
where $A$, $B$, $C$ are matrices of suitable size and with at least one
entry.
Since $w$ is
unitary, $\bar w=w^c$ or $S(w) = w^* $ or equivalently
$$\pmatrix{S(A)&S(B)\cr 0&S(C)} = \pmatrix{A^*&0\cr B^*&C^* }.$$
Therefore $B=0$, moreover $A$ and $C$ are unitary and~$w$ a direct sum
of them (notice that $L^\perp$ is also invariant and $C=w\vert_{L^\perp}$).
By induction hypothesis, the \cors $A$ and $C$ of
dimensions less than~$d$ are direct sums of
irreducible unitary \cors, whence~$w$ is a direct sum of irreducible unitary
\cors.
\item (c) This follows from (a),  (b), and \Ref{P2.15}, part (b).
\item (d) and (e)  follow from \Ref{P2.15}, part (c).
\Ende
Remark: All irreducible \cors can be obtained by decomposition
of tensor products of  those unitary \cors which generate~$\cal A$
as algebra (cf. \Ref{P6.2}, part (a)).

\Abschnitt{Peter-Weyl Theory and Haar measure}
Let $\cal A$ be a compact Hopf $*$-algebra.
Let $\cal I$ be an index set and let $\{\ua..\mid \alpha\in \cal I\}$
be a complete set of non-equivalent irreducible unitary \cors.
Let $I:=u^0$ be the one dimensional \cor.
Then the elements $\ua mn$ form a basis of $\cal A$ (\Ref{P6.2}, part (c)).
\Def{P6.5o} The {\sl Haar measure} is a linear functional on~$\cal A$
defined by
$$h(\ua mn)=\delta_{\alpha,0}.$$
Since the $\ua mn$ are matrix elements of \cors, for all $x\in \cal A$
the Haar measure satisfies the equations
$$(h\tensor \id_{\cal A})\Delta(x)=(\id_{\cal A}\tensor h)\Delta(x)=
h(x)1, h(1)=1. \Eqno{G65}$$
(By definition, also $h(S(x))=h(x)$ holds for all $x\in \cal A$.)
\Ende

In order to compute~$h$ on products, some
preparation is necessary.

\Lemma{P6.5p} For each $\alpha\in \cal I$ there is a strictly positive
definite matrix~$F_\alpha$ such that $(\ua.. )^{cc} = F_\alpha \ua.. F_\alpha\me$.
\Bew
For each $\alpha \in\cal I$, the matrix $\overline{\ua..}$ is also
a \cor and equivalent to a unitary one, say $\ub..$:
$Q_\alpha \overline{\ua..}Q_\alpha\me = \ub..$.
Then $\overline{\ub..} = (\ub.. )^c$ and
$$(\ua.. )^{cc} = (\overline{\ua..})^c =
(Q_\alpha\me \ub.. Q_\alpha)^c = Q_\alpha^T (\ub.. )^c (Q_\alpha\me)^T =
Q_\alpha^T \overline{\ub..} (Q_\alpha\me)^T =
Q_\alpha^T \overline{Q_\alpha} \ua.. \overline{Q_\alpha\me}
(Q_\alpha\me)^T$$
and therefore $(\ua.. )^{cc} = F_\alpha \ua.. F_\alpha\me$ where
$F_\alpha=Q_\alpha^T(Q_\alpha^T)^*$ is a strictly positive definite matrix.
\Ende

Fix an irreducible \cor $v$ and let $n:=\dim(v)$.
Since $S(v)$ is the inverse of $v$, there are intertwiners
$$AI=(v\tensor v^c)A, \ B(v^c\tensor v)=IB,$$
where $A=\sum_{k=1}^n e_k\tensor e_k$ and $B=\sum_{k=1}^n e_k'\tensor e_k'$.

\Lemma{P6.5} Let $v$ and $w$ be irreducible representations of dimensions
$n$ and~$m$ respectively. Then
\item (a) $\Mor(v^c\tensor w, I)\cong \Mor(w,v)$, \
$\Mor( v^c\tensor v, I) = \C B$.
\item (b) $\Mor(I, w\tensor v^c)\cong \Mor(v,w)$, \
$\Mor( I, v\tensor v^c) = \C A$.
\Bew \voritem (a) If $X$ intertwines $v^c\tensor w$ with $I$ then
$X(v^c\tensor w)=IX$ and
$$(\Eins_n\tensor X)(A\tensor\Eins_m)(I\tensor w) =
(\Eins_n\tensor X)(v\tensor v^c \tensor w)(A\tensor \Eins_m)=
(v\tensor I)(\Eins_n\tensor X)(A\tensor \Eins_m).$$
Since $I\tensor w\cong w$ and $v\tensor I\cong v$,
$(\Eins_n\tensor X)(A\tensor\Eins_n)$ can be regarded as intertwiner
of~$w$ and~$v$.
Conversely, let $Y\in \Mor(w,v)$.
Then $Yw=vY$ and
$$B(\Eins_n\tensor Y)(v^c\tensor w) = B(v^c\tensor v)(\Eins_n\tensor Y)=
IB(\Eins_n\tensor Y).$$
Therefore $B(\Eins_n\tensor Y)$ intertwines $v^c\tensor w$ with~$I$.
The maps between $\Mor(v^c\tensor w, I)$ and $\Mor(w,v)$
are inverses of each other because
$(\Eins_n\tensor B)(A\tensor \Eins_n)=
(B\tensor \Eins_n)(\Eins_n\tensor A)=\Eins_n$.
The second statement follows from the first with Schur's \Ref{P2.14}.

\item (b) is proved in a similar way.
\Ende

Now the Haar measure is computed on certain products of basis elements:
\Thm{P6.5a} The Haar measure satisfies the {\sl Peter-Weyl-Woronowicz relations}:
$$h(\ua mn\ub jl^* ) = \delta_{\alpha, \beta}
{ (F_\alpha)_{ln} \delta_{mj} \over Tr(F_\alpha) } \Eqno{G65a}$$
and
$$h(\ub jl^*\ua mn) = \delta_{\alpha, \beta}
{ (F_\alpha\me)_{mj}\delta_{ln}\over Tr(F_\alpha\me)} \Eqno{G65b}$$
for all $\alpha, \beta\in \cal I$, $1\le m,n\le \dim(\ua.. )$,
$1\le j,l\le \dim(\ub.. )$.
\Bew
Let $w$ be any corepresentation (of dimension~$N$).
Application of $h\tensor \id$ and $\id\tensor h$ to
$\Delta w_{ij}$ yields together with \Refg{G65}
$$h(w)w = wh(w) = h(w)1.$$
This matrix equation means
$$\sum_{k=1}^N h(w)_{ik}w_{kj} = h(w)_{ij}1 = \sum_{k=1}^N w_{ik}h(w)_{kj}$$
or equivalently that for each~$i$ the row vector with
coordinates $h(w)_{ij}$ for $j=1, \ldots, N$ intertwines~$w$ with~$I$ and
for each~$j$ the column vector with coordinates $h(w)_{ij}$ for
$i=1, \ldots, N$ intertwines~$I$ with~$w$. These facts will be applied
to the sets $\Mor(I, \uab . )$ and $\Mor(\uab c, I)$ for $\alpha, \beta\in \cal I$.

Therefore for $w=\uab.$ and for fixed indices $k$, $l$, the element
$(h(\ua ik\ub jl^c))_{1\le i,j\le N}$ is in  $\Mor(I, \uab. )$.
By \Ref{P6.5} it vanishes for $\alpha\ne \beta$ and is a multiple
of~$A$ for $\alpha=\beta$. Thus there are numbers $\lambda_{kl}^\alpha\in \C$
such that
$$h(\ua ik\ub jl^c) = \delta_{\alpha, \beta}\lambda_{kl}^\alpha\delta_{ij}\Eqno{G6.1}$$
for all $i$, $j$, $k$, $l$.
Similarly, for $w=\uab c$ and for fixed indices $i$, $j$, the element
$(h(\ua ik^{cc}\ua jl^c))_{1\le k,l\le N}$ is in
$\Mor(\uab c, I)=\C B$. Therefore there
are numbers $\rho_{ij}^\alpha\in\C$ such that
$$h(\ua ik^{cc}\ua jl^c) =
\rho_{ij}^\alpha\delta_{kl}.\Eqno{G6.2}$$
But from \Ref{P6.5p},
$\ua.. = F_\alpha\me (\ua.. )^{cc} F_\alpha$, which yields by
linearity and \Refg{G6.2} the equation
$$h(\ua mn\ua jl^c) = \sum_{i,k} (F_\alpha\me)_{mi}
h(\ua ik^{cc}\ua jl^c) (F_\alpha)_{kn} =
(F_\alpha)_{ln}
\sum_{i} (F_\alpha\me)_{mi} \rho_{ij}^\alpha.$$
Comparison with \Refg{G6.1} and $\ub jl^c = (\ub jl)^*$ yields
$$h(\ua mn\ub jl ^* ) = c_\alpha\delta_{\alpha, \beta}(F_\alpha)_{ln}\delta_{mj}$$
for some $c_\alpha \in\C$. These constants can be evaluated using
the unitarity of~$\ua..$:
$$1=h(1)=\sum_{n} h(\ua mn\ua mn^* ) = c_\alpha\sum_{n} (F_\alpha)_{nn} =
c_\alpha Tr(F_\alpha).$$
The trace of $F_\alpha$ is positive because $F_\alpha$ is positive
definite. This proves \Refg{G65a}. The other equation is proved
in a similar way.
\Ende
\Rema{P6.6} \voritem (a) Since the matrices $F_\alpha$ can be scaled by
a positive number, we {\sl normalize} them
by the condition $Tr(F_\alpha) = Tr(F_\alpha\me)$. After normalization
they are uniquely determined.
\item (b) Example: In the
standard deformation
$SU_q(2)$ for $q\in \R\setminus\{0\}$,
$$F_0 = (1), \ F_{1/2} = \pmatrix{\vert q\vert\me&0\cr 0&\vert q\vert}.$$

\Bew
%\voritem (a) In the following let $m$, $n$ always denote non-negative
%half integers.
%For $m=\halb$, $w=w^m$ is the fundamental corepresentation
$w^{1/2} = w=\Smatrix{\alpha\atop\gamma\cr\beta\atop\delta}$, and
$S(w) = \Smatrix{\delta\atop -q\gamma\cr -\qm\beta\atop\alpha}$.
Then
$$w^{cc} = S^2(w) = \pmatrix{\alpha&q^{-2}\beta\cr q^2\gamma&\delta}
=F_{1/2} w F_{1/2}\me$$
where $F_{1/2}$ is as desired. Note that the absolute value of $q$ must be used,
because the eigenvalues of a positive definite matrix must be positive.
\Ende

\doppelpunktfalse\Thm{P6.7} (Positivity of the Haar measure)

For all $x\in \cal A$, $h(x^*x)\ge 0$, and equality only holds for~$x=0$.
\Bew Since $\cal A$ has a basis $\{ \ua mn \mid 1\le m,n\le \dim(\ua.. ), \alpha\in I\}$,
a general element~$a$ of $\cal A$ can be written as
$$a=\sum_{m,n, \alpha} a_{mn}^\alpha \ua mn.$$
By the second Peter-Weyl-Woronowicz relation (\SRefg{G65b})
$$h(a^* a) = \sum_{\alpha, m,n,p} {(a_{mp}^\alpha(F_\alpha\me)_{mn}
\bar a_{np}^\alpha ) \over Tr(F_\alpha) },$$
in which the sums $\sum_{ m,n} a_{mp}^\alpha(F_\alpha\me)_{mn}
\bar a_{np}^\alpha$ are strictly positive unless all
coefficients~$a_{mp}^\alpha$ for fixed $\alpha$, $p$
vanish, because the matrices~$F_\alpha\me$
are strictly positive definite for all $\alpha$.
\Ende

\doppelpunktfalse\Cor{P6.6n} (Scalar product)

There is a scalar product on $\cal A$ defined by
$(a\mid b):= h(a^*b)$ for all $a,b\in \cal A$.
\Bew This inner product is antilinear in the first argument and linear
in the second argument by definition and positive definite by
\Ref{P6.7}.
\Ende

\doppelpunktfalse\Cor{P6.6o} (Modular Homomorphism)

There is a uniquely determined algebra automorphism $\sigma$ of $\cal A$ such that
$h(ab)=h(b\sigma(a))$ for all $a,b\in\cal A$. It is defined on
elements of the basis as
$$\sigma(\ua mn)= (F_\alpha \ua .. F_\alpha)_{mn}.$$
\Bew \voritem {} Uniqueness: Let $a$ be an element of $\cal A$ and let $a', a''\in\cal A$
such that for all $b\in \cal A$ the equation
$$h(ab)=h(ba') = h(ba'')$$
holds. Then $h(b(a'-a''))=0$ for all $b\in \cal A$, whence $a'=a''$ by
\Ref{P6.6n}.
\item {} Existence:
{}From the second Peter-Weyl-Woronowicz relation it follows that
$$h(\ua jl^* \sigma(\ub mn )) = { \delta_{\alpha, \beta}
(F_\alpha)_{ln}\delta_{mj} \over Tr(F_\alpha)} =
h(\ub mn \ua jl^* ).$$
Therefore by linearity $h(ab)=h(b\sigma(a))$ for all $a, b\in \cal A$.
Moreover $F_0=(1)$ implies $\sigma(1)=1$, and
for all $a,b,c\in \cal A$,
$$h(a\sigma(bc))=h(bca)=h(ca\sigma(b))=h(a\sigma(b)\sigma(c)).$$
Therefore $\sigma$ is a unital algebra homomorphism. Since $F_\alpha$ is
invertible, also $\sigma$ is invertible.
\Ende

\Abschnitt{$C^*$-structure}
For any Hilbert space~$H$ let $( .\vert. )_H$ denote the inner product
and $B(H)$ the set of bounded linear
operators on~$H$. Then~$B(H)$ is a $*$-algebra.
Let $\cal A$ be a compact Hopf $*$-algebra and consider the set
$$\Pi:=\{ \pi\colon {\cal A}\to B(H) \mid
H \hbox{ Hilbert space, $\pi$ unital $*$-homomorphism}\}$$
(it is enough to consider some
fixed $H$ with  $\dim(H)\ge \dim({\cal A})$ as cardinal numbers, 
thus~$\Pi$ is actually a
set).

Fix $\pi\in\Pi$ and let $H$ be the corresponding Hilbert space.
Let $\ua..$ be a unitary corepresentation of $\cal A$. Then~$\pi(\ua.])$ is
a unitary matrix in~$M_{\dim \ua..}(B(H))$ and
$$\sum_{m} \pi(\ua mn)^*\pi(\ua mn) = 1$$
for all $n\le \dim \ua..$. Therefore for all $x\in H$ and $k\le \dim\ua..$
$$\eqalign{(x\mid x)_H &= \sum_m (\pi(\ua mn)^*\pi(\ua mn) x\mid x)_H =\cr
&=
\sum_m (\pi(\ua mn)x\mid \pi(\ua mn) x)_H \ge (\pi(\ua kn)x\mid
\pi(\ua kn) x),}$$
whence the operator norm $\Vert\pi(\ua kn)\Vert$ is at most~1,
and for each $a=\sum_{\alpha, m,n} a_{mn}^\alpha\ua mn\in {\cal A}$
there is the inequality
$$\Vert\pi(a)\Vert \le \sum_{\alpha, m,n} \vert a_{mn}^\alpha\vert
< \infty.$$
Therefore the following definition is possible:
\DefLemma{P7.1} There is a norm $\normc .$ on $\cal A$ such that for all $a\in\cal A$,
$$\normc a = \sup_{\pi\in \Pi} \Vert \pi(a)\Vert.$$
Moreover this norm satisfies the equations
$\normc {ab} \le \normc a\normc b$, $\normc {a^*} = \normc a$,
$\normc {a^*a}= \normc a^2$
for all $a,b\in \cal A$.
\Bew The main problem is to show $\normc a = 0 \Folgt a=0$ for $a\in\cal A$.
The inner product $( .\mid . )$ on~$\cal A$
induces a norm $\normi .$ (cf. \Ref{P6.6n}).
For each $a\in \cal A$ let~$\pi_0(a)$ denote the operator of
left multiplication by~$x$ on~$\cal A$.
Then for all $x\in\cal A$
$$\sum_{m} \normi {\pi_0(\ua mn)(x)}^2 = h(x^*\underbrace{(\sum_{m} (\ua mn)^*\ua mn)}_{=1}x) =
h(x^*x) = \normi x^2,$$
whence the operator norm $\normi .'$ of $\pi_0(\ua mn)$ is at most~1. For all
$a= \sum_{\alpha, m,n} a_{mn}^\alpha\ua mn\in \cal A$
$$\normi {\pi_0(a)}' \le \sum_{\alpha, m,n} \vert a_{mn}^\alpha\vert.$$
Therefore for each $a\in\cal A$ the operator $\pi_0(a)$ is bounded on
$\cal A$ and can be extended to the completion~$H$ of $\cal A$
with respect to the norm $\normi .$ as a bounded linear operator~$\bar \pi_0(a)$
with same operator norm $\normi {\bar\pi_0(a)}':= \normi {\pi_0(a)}'$.
Therefore $\bar \pi_0\in\Pi$, and
$$\normc a = 0 \Folgt \normi {\pi_0(a)}' = 0 \Folgt
\normi {\pi_0(a)1} = 0 \Folgt \normi a = 0 \Folgt a=0.$$
The other properties of this norm follow from the corresponding properties
of the operator norms of the representations in~$\Pi$.
\Ende
\Def{P7.3} Let $A$ be the closure of $\cal A$ with respect to the
norm $\normc .$. Then $A$ is a $C^*$-algebra by \Ref{P7.1}.
\Ende
The following properties of $C^*$-algebras are useful:
\Prop{P7.4} Let $A$ be a $C^*$-algebra. Then
\item (a) There is a Hilbert space $H$ such that $A$ can be embedded
as closed $*$-subalgebra into~$B(H)$ \Dix{2.6.1}.
\item (b) Let $B$ be another $C^*$-algebra. Then each $*$-homomorphism
from $A$ to $B$ is continuous \Dix{1.3.7}.
\Ende
The comultiplication of $\cal A$ can be extended to a $*$-homomorphism
from~$A$ to $A\ttensor A$, where $A\ttensor A$ denotes the
(topological) tensor product of $C^*$-algebras, defined as follows:
Let $H$ be a Hilbert space and let $\iota\colon A\to B(H)$ be an
embedding of $C^*$-algebras. Then $A\ttensor A$ is identified
with the closure of $(\iota\tensor\iota) (A\tensor A)$
in $B(H\ttensor H)$, where $H\ttensor H$ is the (topological)
tensor product of Hilbert spaces. 
The $C^*$-algebra $A\ttensor A$ does not depend (up to isomorphisms) 
on the embedding~$\iota$
\Dix{2.12.15}.
The map
$${\cal A} \buildrel \Delta\over \longrightarrow
{\cal A}\tensor {\cal A}\hookrightarrow B(H\ttensor H)$$
is a $*$-homomorphism called $\pi_1$. Since $H\ttensor H$ is
a Hilbert space, $\pi_1$ belongs to~$\Pi$ and can be extended
to a $*$-homomorphism on~$A$. It is again called $\Delta$. 

\Def{P7.4a} A {\sl \cmqg\empty} is a pair $(A, \Delta)$ or
shortly $A$ where
\item (a) $A$ is a unital $C^*$-algebra generated by some elements
$u_{ij}\in A$ for $1\le i,j\le N$ and some positive integer~$N$,
\item (b) $\Delta\colon A\to A\ttensor A$ is a
unital $*$-homomorphism such that $\Delta(u_{ij})=
\sum_{k=1}^N u_{ik}\tensor u_{kj}$ for all $i$, $j$,
\item (c) the matrices $u$ and $\bar u$ are invertible.
\Ende
\Rema{P7.5} \voritem (a) Let $\cal A$ be a Hopf $*$-algebra generated
as unital algebra by matrix elements of {\sl one} unitary \cor~$u$
or (equivalently) generated as unital
\hbox{$*$-algebra} by matrix elements of a \cor~$v$ such that
$v$ and~$\bar v$ are equivalent to unitary \cors. Then the
$C^*$-algebra constructed as above is a \cmqg.
\item (b) For all positive integers~$N$ the compact Hopf $*$-algebra
of~$SU_q(N)$ gives rise to a \cmqg.
\item (c) The general example of a \cmqg comes from $C^*$-algebras~$A$ as
in (a) after dividing
by closed two-sided ideals $I\subset \{x\in A\colon h(x^*x)=0\}$ such
that~$\Delta$ induces a $*$-homomorphism $A/I\to A/I\ttensor A/I$.
\Ende
\Thm{P7.6}  Let $A$ be a \cmqg constructed as in
\Ref{P7.5}, part (c).
\item (a) Then $\vert h(x)\vert \le \normc x$ for
all $x\in \cal A$, therefore~$h$ can be extended to a (positive)
continuous functional on~$A$, which will be denoted by~$h$ again.
\item (b) The algebra $\cal A$ is embedded into $A$ (because for all
$x\in {\cal A}\setminus\{0\}$
the inequality ${h(x^*x)>0} $ holds).
\item (c) Any \cor of~$A$ (in the sense $\Delta v_{ab}=\sum_c v_{ac}\tensor
v_{cb}$, $v\me$ exists) has matrix elements in~$\cal A$ and thus~$\cal A$
can be recovered from~$A$ as the span of matrix elements of \cors.
\Ende
\Rema{P7.6y} For $I_1:=\{x\in A\colon h(x^*x)=0\}$
(it is a closed two-sided ideal due to \Wi{p. 656}), $h$ is
faithful on $A/{I_1}$ (i.~e.~$h(x^*x)=0\Folgt x=0$), 
while for $I_2:=\{0\}$, $\eps$ is continuous on $A/I_2\cong A$.
In the case of  $SU_q(2)$, $I_1$ and~$I_2$ coincide, cf. \Pdii{Remark 6}.
\Ende

The notion of \cmqgs generalizes that of algebras of continuous functions
on compact groups of matrices. To be more precise:
Let $G$ be a compact group of matrices. Then there is a Haar measure~$\mu$
on~$G$. There is an inner product on $C(G)$ given by
$$(\chi, \psi) :=\int_G \bar \chi \psi d\mu.$$
for $\chi, \psi\in C(G)$.
The algebra $\Poly(G)$ as in \Ref{P2.3} is a compact Hopf $*$-algebra
(cf. proof of \Ref{P2.4}). The inner product as above can also be
expressed as $h(\chi^*\psi)$. Therefore the completion of
$\Poly(G)$ with respect to the norm $\normi .$ is the same as
$L^2(G)$, and the completion of $\Poly(G)$ with respect to the
norm $\normc .$ is the same as $C(G)$.
Here the comultiplication $\Delta\colon C(G)\to C(G)\ttensor C(G)\cong C(G\times G)$
is given by $\Delta(\chi)(g,h) = \chi(gh)$ for all $g,h\in G$ and
$\chi\in C(G)$ (cf. Chapter \Srefk{P2}).
In the following, each compact topological space is by definition a
Hausdorff space.
There are one-to-one correspondences induced by Gel'fand's theorem:
$$\eqalign{\hbox{compact topological spaces~$X$} &\longleftrightarrow
\hbox{unital commutative $C^*$-algebras $C(X)$}\cr
\hbox{continuous mappings $\lambda\colon X\to Y$} &\longleftrightarrow
\hbox{unital $*$-homomorphisms $\lambda^*\colon C(Y)\to C(X)$}\cr
\hbox{cartesian product $X\times Y$} &\longleftrightarrow
\hbox{topological tensor product $C(X)\ttensor C(Y)$}\cr
\hbox{compact group of matrices $G$} &\longleftrightarrow
\vtop{\hbox{\cmqg $C(G)$}\hbox{for commutative ${\cal A}=\Poly(G)$}}}$$

\subsection{P8}{Actions on Quantum Spaces}
\Abschnitt{Definition and spectral decomposition \Pdii{Section 1}}
This chapter deals with a topological counterpart of right 
comodule algebras.
Let $V$ be a topological vector space and $Z\subset V$ a subset. Then
$\langle Z\rangle$ denotes the closure of the linear span of the
elements of~$Z$ in~$V$.
\Def{P8.1} Let $(A, \Delta)$ be a \cmqg and $B$ a unital $C^*$-algebra.
 The unital $*$-homomorphism $\Gamma\colon B\to B\ttensor A$
 is called a {\sl coaction} for~$A$ on~$B$ if
\item (a) $(\Gamma\tensor \id_A)\Gamma = (\id_B\tensor \Delta)\Gamma$,
\item (b) $B\tensor A
=\langle (\id_B\tensor y)\Gamma(x)\mid x\in B, y\in A\rangle$.
\Ende

\Rema{P8.2} \voritem (a)
Let~$G$ be a compact group of matrices, $X$ a compact topological
space and $X\times G\to X$, $(x,g)\mapsto xg$ for $x\in X$ and
$g\in G$, an action. Then there is a
coaction $\Gamma\colon C(X)\to C(X\times G)$ given by
$\Gamma(\chi)(x,g)=\chi(xg)$ for all $\chi\in C(X)$, $g\in G$, $x\in X$.
The properties
$x(gh)=(xg)h$ and $xe=x$ for all $g,h\in G$, $x\in X$ correspond
to Conditions (a) and (b) in \Ref{P8.1} respectively. Given a coaction
as in \Ref{P8.1} for commutative $A$ and~$B$, the group action can
be recovered by Gel'fand's theorem.
%\item (b) The second condition in \Ref{P8.1} is not necessary if~$A$
%is a Hopf algebra.
\item (b) Quantum analogues of left actions are considered in 
\Pdii{Remark 7}.
\Ende

\Thm{P8.3} Let $A$ be a \cmqg, $B$  a unital $C^*$-algebra and
$\Gamma$ a coaction. Then there exists a maximal $*$-subalgebra $\cal B$
of $B$ such that $\cal B$ is dense in~$B$ and an $\cal A$ right
comodule algebra, i.~e.~for $\gamma:=\Gamma\vert_{\cal B}$:
$$\gamma({\cal B})\subset
{\cal B\tensor \cal A}, \ (\gamma\tensor \id)\gamma = (\id\tensor \Delta)\gamma, \
(\id\tensor \eps)\gamma = \id.$$
For each $\alpha\in \cal I$ there is a set~$I_\alpha$ such that
the algebra~$\cal B$ has a basis $e_{\alpha rk}$ for $\alpha\in \cal I$,
$r\in I_\alpha$, $1\le k\le \dim(\ua.. )$ such that
$$\Gamma(e_{\alpha rk}) = \sum_{s} e_{\alpha rs}\tensor \ua sk.$$
\smallskip
\rm Idea of proof (cf. \Pdii{Theorem 1.5}): From the Peter-Weyl-Woronowicz
relation (\SRefg{G65b}) it follows that there are elements
$x^\alpha_{sm}\in \cal A$ which span~$\cal A$
such that the continuous linear functionals
$$\rho_{sm}^\alpha\colon A\to \C, \ x\mapsto h(x_{sm}^\alpha x)$$
satisfy $\rho_{sm}^\alpha(\ub kr) = \delta_{\alpha, \beta}
\delta_{sk}\delta_{mr}$. Then the
%all $\alpha$, $s$, $m$) span~$\cal A$ and the
operators
$$E_{sm}^\alpha = (\id_B\tensor \rho_{sm}^\alpha)\Gamma\colon B\to B$$
have properties of matrix units. The traces $\sum_s E_{ss}^\alpha$
are projections onto  subspaces~${W_\alpha\subseteq B}$ which contain
all elements~$x\in B$ such that
$$\Gamma(x)\subset B\tensor (\oplus_{ik} \C\ua ik).$$
Construction of the basis: For each $\alpha\in \cal I$ let $\{e_{\alpha r1}\mid
r\in I_\alpha\}$ be a basis of the vector space $\Im(E_{11}^\alpha)$
and $e_{\alpha rs}:= E_{s1}^\alpha(e_{\alpha r 1})$. Let $\cal B$ denote the
linear span of all elements $e_{\alpha rs}$.
Then the closure of~$\cal B$ is
$$\eqalign{\langle E_{sm}^\alpha(x)\mid x\in {\cal B}, \alpha, s,m \rangle &=
\langle E_{sm}^\alpha(x)\mid x\in B, \alpha, s, m\rangle = \cr
&=\langle (\id\tensor h)(\id\tensor x_{sm}^\alpha)\Gamma(x) \mid
x\in B, \alpha, s,m\rangle =\cr
&=\langle (\id\tensor h)\langle (\id\tensor y)\Gamma(x)\mid y\in A,
x\in B\rangle\rangle =\cr
&= \langle(\id\tensor h)(B\ttensor A)\rangle = B.}$$
\Ende

\Def{P8.3a} Let a \cmqg $A$ coact by~$\Gamma$ on a quantum space~$B$.
\item (a)
For each $\alpha\in \cal I$, the number $c_\alpha$ denotes the
cardinality of $I_\alpha$ as in \Ref{P8.3} and is called \Quote{multiplicity
of~$\ua ..$ in the spectrum of~$\Gamma$}.
\item (b) For each $\alpha\in \cal I$ let $W_\alpha$ be the linear
span of the elements $e_{\alpha rs}$ as in \Ref{P8.3}.
\Ende

\Abschnitt{Quantum spheres \Pd{}}
Since the quantum groups $SU_q(2)$ and $SU_{1/q}(2)$ are isomorphic by
\Ref{P3.2}, we can restrict ourselves to the case $q\in [-1, 1]\setminus\{0\}$.
For the quantum $SU(2)$ groups, $\cal I$ is the set of non-negative
half integers and $u^k=w^k$ for $k\in\cal I$.
We want to classify coactions~$\Gamma$ of $SU_q(2)$ such that
\item (1) $\displaystyle c_k = \cases{ 1& if $k\in \N_0$\cr
0& if $k\in \N_0+\halb,$}$
\item (2) the subspaces $W_0$ and $W_1$ generate $B$ as a $C^*$-algebra.

The pairs $(B, \Gamma)$ are called \Quote{quantum spheres}
(cf. the case $q=1$ in \Ref{P8.9} below).
For convenience, the matrix elements of the unitary irreducible
\cors of $SU_q(2)$ will be indexed by numbers in the index set
$$N_\alpha:=\{-\alpha, -\alpha+1, \ldots , \alpha\}$$
instead of the index set $\{1, \ldots, 2\alpha+1\}$ for
each $\alpha\in \halb \N_0$.

\doppelpunktfalse\Thm{P8.9} {\rm \Pd{}:}
In the case $q=1$ there is only one object $B=C(S^2)$ and the coaction
is induced by the standard right action of $SU(2)$ on the sphere~$S^2$.
Here $W_0=\C\, 1$ and $W_1=\C x+\C y+\C z$. Then Condition~(2)
means that the coordinates $x$, $y$, $z$ separate the points of~$S^2$
by the Stone-Weierstrass theorem.

In the case $q=-1$ there is only one object $B_{-1,0}$ with coaction
$\Gamma_{-1,0}$.

In the case $-1<q<1$ and $q\ne 0$ there are---up to isomorphisms---the following quantum
spaces~$B_{qc}$ for $c\in \R_0^+\cup \{\infty\}$.
The $C^*$-algebra $B_{qc}$ is
generated by the elements $e_{-1}$, $e_0$, $e_1$ of $W_1$ subject to
the relations
$$e_i^* = e_{-i} \hbox{ for $i\in \{-1,0,1\}$,}$$
$$\left.\eqalign{
(1+q^2)(\em e_1+q^{-2}e_1\em)+e_0^2 &= \rho 1,\cr
e_0\em -q^2\em e_0 &= \lambda\em \cr
(1+q^2)(\em e_1 - e_1\em)+(1-q^2)e_0^2 &= \lambda e_0,\cr
e_1e_0-q^2 e_0e_1 &= \lambda e_1,}\right\}\Eqno{G16}$$
where
$$\lambda = \cases{1-q^2& if $c\in \R$\cr 0& if $c=\infty$} \quad\hbox{ and }
\rho = \cases{(1+q^2)^2q^{-2}c+1& if $c\in \R$\cr
(1+q^2)^2q^{-2}& if $c=\infty$.}$$
The coaction $\Gamma_{qc}$ is given by
$$\Gamma(e_i) = \sum_{j=-1}^1 e_j\tensor u_{ji}^1$$
for $i\in \{-1,0,1\}$. Here we choose a non-unitary form
$$u^1=\pmatrix{\delta^2&-(q^2+1)\delta\gamma&-q\gamma^2\cr
-\qm\beta\delta& 1+(q+\qm)\beta\gamma&\alpha\gamma\cr
-\qm\beta^2&(q+\qm)\beta\alpha&\alpha^2\cr}.$$
\smallskip\rm
Ideas of proof: Due to \Ref{P8.3} and Condition~(1), the algebra $\cal B$ has the linear
basis ${\{e_{\alpha k}\mid \alpha\in \N_0, k\in N_\alpha\}}$
such that
$$\Gamma(e_{\alpha k}) = \sum_{s\in N_\alpha} e_{\alpha s}\tensor
\ua sk \hbox{ for $\alpha\in \N_0$, $k\in N_\alpha$}.$$
Therefore the $e_{\alpha k}$'s are analogues of spherical harmonics.
One has $(u^1_{lk})^* = u^1_{-l,-k}$.
Then
$$\Gamma(e_{-k}^* )= \sum_l e_{-l}^*\tensor (u^1_{-l,-k})^* =
\sum_l e_{-l}^* \tensor u^1_{lk}.$$
{}From the irreducibility of $u^1$ it follows that there is a constant~$c$ such
that $e_{-k}^*= ce_k$ for all~$k$. Moreover the modulus of~$c$ is one because
of $e_k=(e_k^* )^* = (ce_{-k})^* = c\bar c e_k$. Thus it is possible to
achieve $c=1$ by scaling the elements~$e_k$ with a suitable complex number
of modulus one.

Now consider products of the generators: 
Because of the Clebsch-Gordan relation $u^1\tensor u^1\cong u^0\oplus u^1\oplus
u^2$ there are injective intertwiners~$G^\alpha\in \Mor(u^\alpha, u^1\tensor u^1)$ 
for $\alpha\in\{0,1,2\}$. From the equation
$$\Gamma(e_ke_l)= \sum_{m,r} e_me_r\tensor u^1_{mk}u^1_{rl}$$
it follows for the elements $\tilde e_{\alpha,t}:=\sum_{k,l}e_ke_l G^\alpha_{kl,t}$:
$$\Gamma(\tilde e_{\alpha,t})=\sum_{k,l,m,r} e_me_r\tensor
u^1_{mk}u^1_{rl}G^\alpha_{kl,t} =
\sum_n \Bigl(\underbrace{\sum_{m,r} e_me_rG^\alpha_{mr,n}}_{\textstyle=\tilde e_{\alpha,n}}\Bigr)\tensor u^\alpha_{nt}.$$
Therefore the elements $\tilde e_{\alpha,t}$ satisfy the same relations for
the coaction as the elements~$e_k$. Since the \cors $u^\alpha$ are irreducible,
there are constants $\lambda_\alpha\in\C$ such that
$\tilde e_{\alpha,t}=\lambda_\alpha e_{\alpha,t}$.
For $\alpha\in\{0,1\}$ this gives relations for the
generators:
$$\eqalign{\sum_{k,l} e_ke_lG^1_{rl,t}&=\lambda e_t \hbox{ (here $\lambda=\lambda_1$),}\cr
\sum_{k,l} e_ke_l G^0_{rl,0}&= \rho 1\hbox{ (here $\rho=\lambda_0$).}\cr}$$
These are the relations (\SRefg{G16})
for the quantum spheres. Applying \Quote{$*$} to both sides,
we obtain that $\lambda$ and~$\rho$ are real.
There is still the freedom of scaling the $e_k$'s by a non-zero real number.
Consider the case $0<\vert q\vert<1$.  If $\lambda$ does not vanish, it
can be scaled to the value
$\lambda=1-q^2$. Then define~$c$ by
$$\rho=(1+q^2)^2q^{-2}c+1.$$
The existence of a faithful $C^*$-norm on $\cal B$ implies that $c$ is
a non-negative number.
It remains~$\lambda=0$, $\rho$ positive ($B$ is a $C^*$-algebra).
Then~$\rho$ can be scaled to the value $(1+q^2)^2q^{-2}$.

These $(B, \Gamma)$'s are indeed quantum spheres. No extra relation can be imposed,
because then we would get a coaction for a quantum subspace.
But $c_0=1$ means that the space
is homogeneous (cf. \Pdii{Definition 1.8}), and from the facts that~$h$
is faithful (i.~e.~$h(x^*x)=0 \Folgt x=0$) and the counit is continuous
(cf. \Ref{P7.6y}) it follows here that the homogeneous space
corresponding to~$B$ has no non-trivial homogeneous subspaces
(this idea stands behind the proof in the paper \Pd{}).

The case $q=1$ can be handled similarly, and the case $q=-1$ reduces to~$q=1$.
\Ende

\Rema{P8.10} \voritem (a) If the first condition for the quantum spheres
is weakened to $c_0=c_1=1$, there are some more
homogeneous spaces for $c\in \{c(2), c(3), \ldots\}$, $0<\vert q\vert <1$,
where
$$c(n)=  -q^{2n}/ (1+q^{2n})^2 \hbox{ for all $n\in\N$.}$$
These objects satisfy the conditions
$$c_k=\cases { 1&if $k=0,1,\ldots, n-1$\cr 0&otherwise.}$$
There exist analogues of these objects in the
case $q=1$.
They correspond (cf. \Pd{}) to the adjoint action of $SU(2)$ on
$U(\frak {su}(2))$ taken in its $n$-dimensional irreducible
$*$-representation ($X^*=-X$ for $X\in \frak {su}(2)$).
\item (b) For $0<\vert q\vert<1$, $c\in \R_0^+\cup\{\infty\}\cup
\{c(2), c(3), \ldots\}$ the quantum sphere $S^2_{qc}=(B_{qc}, \Gamma_{qc})$
is a
quotient space if and only if $c=0$, embeddable (i.~e.~can be regarded
as a non-zero $C^*$-subalgebra of~$A$,
where $\Gamma$ is induced by the comultiplication) if 
and only if $c\in [0, \infty]$,
and homogeneous for all considered $c$
(for the compact groups of matrices these three
notions coincide).
\item (c) An algebraic version of \Ref{P8.9} can be found in \Schm{}.
\Ende

\subsection{P9}{Quantum Lorentz groups (cf. \WZ{})}
The algebra ${\cal A}=\Poly(SL(2, \C))$ is called the algebra of polynomials
on the Lorentz group. Its \cors have the following properties (cf. Chapter
\Srefk{P3}):
\item (1) There are irreducible \cors $w^\alpha$ for
$\alpha\in\halb\N_0$ such that all non-equivalent
irreducible \cors are
$w^\alpha\tensor \overline{w^\beta}$ for $\alpha, \beta\in \einhalb \N_0$.
\item (2) $\dim(w^\alpha) = 2\alpha+1$ for all $\alpha$,
\item (3) $w^\alpha\tensor w^\beta \cong w^{\vert \alpha-\beta\vert}\oplus
w^{\vert \alpha-\beta\vert+1}\oplus\cdots w^{\alpha+\beta}$
(Clebsch Gordan),
\item (4) Each \cor is completely reducible, or equivalently,
the matrix elements $w_{ij}^\alpha (w_{kl}^\beta)^*$ give
a basis of $\cal A$.
\item (5) For all $\alpha , \beta\in \halb \N_0$ the \cors
$w^\alpha\tensor \overline{w^\beta}$ and $\overline{w^\beta}\tensor w^\alpha$
are equivalent.
\Def{P9.1o} A {\sl quantum Lorentz group} is a Hopf $*$-algebra~$\cal A$
satisfying
properties (1)--(5).
\Ende

\Thm{P9.1} Up to isomorphisms, all quantum Lorentz groups~$\cal A$
are given as follows:
The Hopf $*$-algebra $\cal A$ is
generated by the matrix elements $w_{ij}$ ($1\le i,j\le 2$)
of the fundamental \cor $w:=w^{1/2}$ and relations
\item(i) $(w\tensor w)E=E$,
\item (ii) $E'(w\tensor w)=E'$,
\item (iii) $X(w\tensor \bar w)=(\bar w\tensor w)X$,

where the base field~\C is canonically embedded into~$\cal A$, the
vectors $E'\in \C^2\tensor\C^2$ and
$E\in {}^2\C\tensor{}^2\C$ are the same as in \Ref{P3.2}
and $X\in M_4(\C)$
satisfies the properties:
\item (iv) $X$ is invertible,
\item (v)  there is a scalar factor $c\in \C\setminus\{0\}$ such that
$\tau \bar X\tau = cX$,
\item (vi) the intertwiners $\Eins_2\tensor E$ and
$(X\tensor \Eins_2)(\Eins_2\tensor X)(E\tensor \Eins_2)$
in $\Mor(\bar w, \bar w\tensor w\tensor w)$ are proportional (note
that $\bar w\cong w^0\tensor \bar w\cong \bar w\tensor w^0$).
\smallskip
\rm Idea of proof: Necessity of relations: Restrict attention to
the \cors $w^\alpha$ first. Their
matrix elements give a basis of a
quantum $SL(2)$-group~$\cal H$ as in \Ref{P3.2}.
This shows conditions~(i) and~(ii) and gives~$E$ and~$E'$.
{}From assertions (1) and (4) it follows that
there is a linear isomorphism
$${\cal A}\cong {\cal H}\cdot{\cal H}^* \cong {\cal H}\tensor {\cal H}^*,
\quad w^\alpha_{kl}(w^\beta_{mn})^*
\mapsto w^\alpha_{kl}\cdot(w^\beta_{mn})^*\mapsto w^\alpha_{kl}\tensor_{\Bbb C}
(w^\beta_{mn})^*,$$
where \Quote{$\cdot$} denotes multiplication.
Assertion (5) for $\alpha=\beta=\halb$ shows that there is a bijective
intertwiner $X\in \Mor(w\tensor \bar w, \bar w\tensor w)$, which gives
conditions~(iii) and~(iv). Apply the map \Quote{$*$} to~(iii) and
use the formula $\overline{v\tensor w}=\tau(\bar w\tensor \bar v)\tau$
as in the proof of \Ref{P5.1o}, part (c):
$$\bar X\overline{(w\tensor \bar w)} = \overline{(\bar w\tensor w)}\bar X \Folgt
\bar X\tau (w\tensor \bar w)\tau= \tau(\bar w\tensor w)\tau\bar X
\Folgt \tau\bar X\tau ( w\tensor \bar w) =
(\bar w\tensor w)\tau\bar X\tau.$$
Since $\bar w\tensor w$ and $w\tensor \bar w$ are irreducible, the
intertwiners $X$ and~$\tau\bar X\tau$ must be proportional, which
gives Condition~(v).
The last condition follows, because both $\Eins_2\tensor \Eins_2$
and
$$X(\Eins_2\tensor\Eins_2\tensor E')(\Eins_2\tensor X\tensor \Eins_2)(E\tensor\Eins_2\tensor \Eins_2)$$
are elements of $\Mor(\bar w\tensor w)$.

Existence: We set ${\cal A}:= {\cal H}\tensor{\cal H}^*$
with~$\cal H$ as in \Ref{P3.2} and laborously introduce the Hopf
$*$-algebra stucture on~$\cal A$ by means of (iii)--(vi).

Sufficiency of relations: More relations would make the elements 
$w^\alpha_{ij}(w^\beta_{kl})^*$
linearly dependent.
\Ende

\Rema{P9.2} \voritem (a) Possible matrices~$X$ have been found (up to isomorphisms of the corresponding Hopf $*$-algebras)
 in \WZ{}.
\item (b) There is also a topological structure for
two examples of~$\cal A$ (\WP{}, \WZi{}) which uses the notion
of affiliated elements \Wv{}.
\item (c) Quantum Poincar\'e groups arise by adding translations \Pdvi{}.
\item (d) Quantum analogues of $\Poly(SL(N, \C))$ were considered
in \Pdiii{} (cf. \Zal{}).
\Ende

\advance\baselineskip by -0.3mm

\subsection{P10}{References}
\bigskip
\Litverz
\bye